\documentclass[12pt]{article}
\begin{document}

\newtheorem{example}{Example}[section]
\newtheorem{theorem}{Theorem}[section]
\newtheorem{lemma}{Lemma}[section]
\newtheorem{corollary}{Corollary}[section]
\newtheorem{proposition}{Proposition}[section]
\newtheorem{remark}{Remark}[section]
\setlength{\textwidth}{13cm}
\def \halmos{\hfill\mbox{ qed}\\}
\newcommand{\eqnsection}
{\renewcommand{\theequation}{\thesection.\arabic{equation}}
\makeatletter \csname  @addtoreset\endcsname{equation}{section}
\makeatother}

%

\def \njc{{\bf !!! note Jay's change  !!!}}
\def  \enc{{\bf end of current changes }}
\def \nnc{{\bf !!! note NEW	 change  !!!}}
\def\square{{\vcenter{\vbox{\hrule height.3pt
                \hbox{\vrule width.3pt height5pt \kern5pt
                   \vrule width.3pt}
                \hrule height.3pt}}}}
\def \grad{\bigtriangledown}
\def \nc{{\bf !!! note change !!! }}
    \def \Proof{\noindent{\bf Proof}$\quad$}
\newcommand{\re}[1]{(\ref{#1})}
\def \ov{\overline}
\def \un{\underline}
\def \be{\begin{equation}}
\def \ee{\end{equation}}
\def \bt{\begin{theorem}}
\def \et{\end{theorem}}
\def \bc{\begin{corollary}}
\def \ec{\end{corollary}}
\def \br{\begin{remark} }
\def \er{ \end{remark}}
\def \bl{\begin{lemma}}
\def \el{\end{lemma}}
\def \bex{\begin{example}}
\def \eex{\end{example}}
\def \bea{\begin{eqnarray}}
\def \eea{\end{eqnarray}}
\def \bas{\begin{eqnarray*}}
\def \eas{\end{eqnarray*}}
\def \al{\alpha}
\def \bb{\beta}
\def \ga{\gamma}
\def \Ga{\Gamma}
\def \de{\delta}
\def \De{\Delta}
\def \ep{\epsilon}
\def \vep{\varepsilon}
\def \la{\lambda}
\def \La{\Lambda}
\def \ka{\kappa}
\def \om{\omega}
\def \Om{\Omega}
\def \va{\varrho}
\def \ffi{\Phi}
\def \vf{\varphi}
\def \si{\sigma}
\def \Si{\Sigma}
\def \vsi{\varsigma}
\def \th{\theta}
\def \Th{\Theta}
\def \ups{\Upsilon}
\def \ze{\zeta}
\def \tr{\nabla}
\def \ff{\infty}
\def \wh{\widehat}
\def \wt{\widetilde}
\def \dar{\downarrow}
\def \rar{\rightarrow}
\def \uar{\uparrow}
\def \sbs{\subseteq}
\def \mpt{\mapsto}
\def \R{{\bf R}}
\def \G{{\bf G}}
\def \H{{\bf H}}
\def \Z{{\bf Z}}
\def \S{{\bf S}}
\def \sfB{{\sf B}}
\def \sfS{{\sf S}}
\def \T{{\bf T}}
\def \C{{\bf C}}
\def \AA{{\mathcal A}}
\def \BB{{\mathcal B}}
\def \CC{{\mathcal C}}
\def \DD{{\mathcal D}}
\def \EE{{\mathcal E}}
\def \FF{{\mathcal F}}
\def \GG{{\mathcal G}}
\def \HH{{\mathcal H}}
\def \II{{\mathcal I}}
\def \JJ{{\mathcal J}}
\def \KK{{\mathcal K}}
\def \LL{{\mathcal L}}
\def \MM{{\mathcal M}}
\def \NN{{\mathcal N}}
\def \OO{{\mathcal O}}
\def \PP{{\mathcal P}}
\def \QQ{{\mathcal Q}}
\def \RR{{\mathcal R}}
\def \SS{{\mathcal S}}
\def \TT{{\mathcal T}}
\def \UU{{\mathcal U}}
\def \VV{{\mathcal V}}
\def \ZZ{{\mathcal Z}}
\def \Pxh{P^{x/h}}
\def \Exh{E^{x/h}}
\def \Px{P^{x}}
\def \Ex{E^{x}}
\def \Prh{P^{\rho/h}}
\def \Erh{E^{\rho/h}}
\def \p{p_{t}(x,y)}
\def \({\left(}
\def \){\right)}
\def \lk{\left[}
\def \rk{\right]}
\def \lc{\left\{}
\def \rc{\right\}}
\def \bsq{\ $\Box$}
\def \nn{\nonumber}
\def \Bo{\bigotimes}
\def \bo{\times}
\def \ot{\times}
\def \c*{ \begin{center}***********************\end{center}}
\def \bs{\begin{slide} }
\def \es{\end{slide} }
\def \bpr{\begin{proof} }
\def \epr{\end{proof} }
\def \cd{\,\cdot\,}
\def \st{\stackrel{def}{=}}
\def \as{almost surely }
\def \fix {{\bf  !!!!! }}
\def \stl{\stackrel{law}{=}}
\def \std{\stackrel{dist}{=}}
\def \stdto{\stackrel{dist}{\longrightarrow}}
\def  \enc{{ \bf end of current changes}}
\def\square{{\vcenter{\vbox{\hrule height.3pt
                \hbox{\vrule width.3pt height5pt \kern5pt
                   \vrule width.3pt}
                \hrule height.3pt}}}}
\def\qed{{\hfill $\square$ }}

       \def \tb{|\!|\!|}

    \eqnsection
\bibliographystyle{amsplain}

\title{Existence of a critical point for the infinite  divisibility of squares  of Gaussian vectors in $R^{2}$ with non--zero mean}
    \author{ Michael B. Marcus\,\, Jay Rosen \thanks{Research of both 
authors supported by  grants from the National Science Foundation and
PSCCUNY.}}

\maketitle

\begin{abstract} Let $G=(G_{1},G_{2})$ be a Gaussian vector in $R^{2}$ with $EG_{1}G_{2}\ne 0$.   Let $c_{1},c_{2}\in R^{1}$. A necessary and sufficient condition for $G=((G_{1}+c_{1}\al)^{2},(G_{2}+c_{2}\al)^{2})$ to be infinitely divisible for all $\al\in R^{1}$ is that 
\be \Ga_{i,i}\ge \frac{c_{i}}{c_{j}}\Ga_{i,j}>0\qquad\forall\,1\le i\ne j\le 2.\label{0.1}
 \ee
In this paper we show that when (\ref{0.1}) does not hold there exists an $0<\al_{0}<\ff $ such that  $G=((G_{1}+c_{1}\al)^{2},(G_{2}+c_{2}\al)^{2})$ is infinitely divisible for all $|\al|\le \al_{0}$ but not for any $|\al|>\al_{0}$.

 \end{abstract}

 \section{Introduction}

Let $\eta=(\eta _1,\ldots,\eta _n)$ be an $R^n$ valued  Gaussian
random variable. $\eta$  is said to have  infinitely  divisible
squares\index{ infinitely divisible squares} if
$\eta^2:=(\eta^{ 2}_1,\ldots,\eta^{ 2}_n)$ is infinitely divisible, i.e. for
any $r$ we can find an
$R^{ n}$ valued random vector $Z_{r }$ such that
\begin{equation}
   \eta^2\stl Z_{ r,1}+\cdots+Z_{ r,r}\label{id.1},
\end{equation}  where $\{Z_{ r,j} \}$,  $j=1,\ldots,r $ are independent
identically distributed copies of
$Z_{ r}$. We   express this by saying that $\eta^2$ is infinitely divisible.

Paul L\'evy proposed the problem of characterizing which Gaussian vectors have infinitely divisible squares. It is easy to see that a single Gaussian random variable
has infinitely divisible squares. However, even for vectors in $R^{2}$ this is a difficult problem. It seems that L\'evy incorrectly conjectured that not all Gaussian vectors in $R^{2}$  have infinitely divisible squares. If he had said $R^{3}$ his conjecture would have been correct.

L\'evy's problem was solved by Griffiths and Bapapt
\cite{Bapat,Griffiths}, (see also \cite[Theorem 13.2.1]{book}).  

\bt\label{lem-Bapat0} Let $G=(G _{  { 1}},\ldots,G _n)$ be a  mean
zero Gaussian random variable with strictly positive definite covariance
       matrix $\Ga=\{\Ga_{i,j}\}= \{E(G _iG _j)\}$.    Then
$G^2$ is infinitely divisible if and only if there exists a signature matrix
      $\mathcal{N}$ such that  
\be
\NN\Ga^{-1}\NN\qquad \quad\mbox{is an
$M$ matrix.} \label{1}
\ee 
\et

We need to define the different types of 
matrices that appear in this theorem.   Let
$A=\{ a_{ i,j}\}_{ 1\leq i,j\leq n}$ be an
$n\times n$ matrix. We call $A$ a positive matrix  and write
$A\geq 0$ if
$a_{ i,j}\geq 0$ for all
$i,j$.  We say that $A$ has positive row sums if $\sum_{j=1}^n
a_{i,j}\ge 0$ for all $1\le i\le n$.

\medskip  The matrix
$A$  is said to be an
$M$ matrix   if
\begin{enumerate}
\item[(1)] $a_{ i,j}\leq 0$ for all $i\neq j$.
\item[(2)] $A$ is nonsingular and $A^{ -1}\geq 0$.
\end{enumerate}

  A  matrix is called a signature matrix if its off--diagonal entries are all zero
and   its diagonal entries are either one or minus one. The role of the signature matrix is easy to understand.  It simply accounts for
the fact that if
$G$ has an infinitely divisible square, then so does
$(\ep_1G_1,\ldots,\ep_nG_n)$
 for any choice of $\ep_i=\pm 1$,
$i=1,\ldots,n$.   Therefore, if (\ref{1}) holds for $\NN$ with diagonal elements
$n_1,\ldots,n_n$
\be
\(\NN\Ga^{-1}\NN\)^{-1}=\NN\Ga \NN  \ge 0\label{2}
\ee
since the inverse of an  $M$ matrix is positive. Thus
$(n_1G_1,\ldots,n_nG_n)$ has a positive covariance   matrix and its inverse is
an
$M$ matrix.  (For this reason, in studying mean zero Gaussian vectors with
infinitely divisible squares one can restrict ones attention to vectors with
positive covariance.)

\medskip	

The natural next step was to   characterize   Gaussian processes with infinitely
divisible squares  which  do  not have mean zero. We set 
$\eta_{i}=G_{ i}+c_{ i}$,  $EG_{ i}=0$, $i=1,,\dots,n$.   Let $\Ga$ be 
the covariance matrix of  $(G _1 ,\dots,  G  _n)$ and set
\be
 c:=(c _1,\dots,c _n ).
\ee
    Set
\be  G+c :=(G _1 +c _1,\dots,G _n +c _n )\label{1.5}
   \ee 
and
\be  (G+c)^2 :=
((G _1 +c _1 )^2,\dots,(G _n +c _n )^2)\label{1.6xx}.
\ee

  Results about  the infinite divisibility of $ (G+c)^2$ when  $c_{1}=\cdots =c_{n}$,  are given in the work of N.
Eisenbaum \cite{Eisen03,Eisen05}   and then in joint work by Eisenbaum 
and H. Kaspi \cite{EK06},  as a by product of their  characterization of Gaussian processes with
a  covariance   that is the 0-potential density of a symmetric Markov process. We point out later in this Introduction how   Gaussian   vectors with infinitely divisible squares are related to the local times of the Markov chain that is determined by the covariance  of the Gaussian vector. It is this connection   between Gaussian vectors with infinitely divisible squares and the local times of Markov chains, and more generally, between Gaussian processes with infinitely divisible squares and the local times of Markov processes, that enhances our interest in the question of characterizing Gaussian vectors with infinitely divisible squares.  

\medskip	 
Some of the results   in  \cite{Eisen03,Eisen05,EK06} are presented and expanded in \cite[Chapter 13]{book}.  The following  theorem  is taken from  
\cite [Theorem 13.3.1 and Lemma 13.3.2]{book}.

\bt\label{theo-book} Let $G=(G _{  { 1}},\ldots,G _n)$ be a  mean
zero Gaussian random variable with strictly positive definite covariance
       matrix $\Ga=\{\Ga_{i,j}\}= \{E(G _iG _j)\}$.  Let ${\bf 1}=(1,\ldots,1)\in R^{n}$.  Then the following are
equivalent:

\begin{itemize}
\item[(1)] $(G+{\bf 1} \al)$ has infinitely divisible squares for all $\al\in R^1$;
\item[(2)] For $\xi=N(0,b^2)$ independent of $G$, $(G_1+ \xi,\ldots, G_n+ \xi,
\xi)$ has infinitely divisible squares   for some $b\neq 0$. Furthermore, if this holds for  some $b\neq 0$, it holds   for all $b\in R^1$, with
$N(0,0)=0$.
\item[(3)] $ \Ga^{-1} $ is an
$M$ matrix with positive row sums.
\end{itemize}
\et

In \cite{IDS}, Theorem \ref{theo-book} is generalized so that the mean of the components of $G+c$ in (\ref{1.5}) need not be the same. In this   generalization certain trivial cases spoil the simplicity of the final result. We avoid them by requiring that the covariance matrix of the Gaussian process is irreducible.

\bt\label{lem-Bapat} Let $G=(G _{  { 1}},\ldots,G _n)$ be a  mean
zero Gaussian random variable with irreducible strictly positive definite covariance
       matrix $\Ga=\{\Ga_{i,j}\}= \{E(G _iG _j)\}$.  Let $c=( c_{
1},
\ldots ,c_{ n})\in R^n$,   $c\ne \mathbf{  0}  $  and let $C$ be a
diagonal matrix with
$c_i=C_{i,i}$, $1\le i\le n$ . Then the following are equivalent: 

\begin{itemize}
\item[(1)] $G+c \al$ has infinitely divisible squares  for all $ \al\in R^1$;
\item[(2)] For $\xi=N(0,b^2)$ independent of $G$, $(G_1+c_1 \xi,\ldots,
G_n+c_n \xi, \xi)$ has infinitely divisible squares   for some $b\neq 0$. Furthermore, if this holds for  some $b\neq 0$, it holds   for all $b\in R^1$;
\item[(3)]   
$C\, \Ga^{-1} \, C$ is an
$M$ matrix with positive row sums.
\end{itemize}
\et

\medskip
By definition, when $(G+c )^2 $ is infinitely divisible, it can
be written as in (\ref{id.1}) as a sum of $r$ independent
identically distributed random variables,  for all $r\ge 1$.   
Based on the work of   Eisenbaum and Kaspi mentioned above and the joint
paper
\cite{fiveauthors} we can actually describe the decomposition. We give a
rough description here. For details see \cite{Eisen03,Eisen05,EK06} and 
\cite[Chapter 13]{book}.

Assume that (1), (2) and (3) of Theorem \ref{lem-Bapat} hold.  Let
\be
{G \over c}=\({G_{1} \over c_{1}},\ldots,{G_{n} \over c_{n}}\).\label{2.1xx}
\ee
 Let
$\Ga_c$ denote the covariance matrix of $G/c$. Theorem
\ref{theo-book} holds for $ G/c$ and
$\Ga_c$, so $\Ga_c^{-1}$ is an $M$ matrix with positive row sums. To
be specific let 
$G/c\in  R^n$. Set $S=\{1,\ldots,n\}$. By \cite[Theorem 13.1.2]{book}, $\Ga_c$
is the 0-potential density of a strongly symmetric   transient Borel right process, say
$X$, on $S$. We show in  the proof of \cite[Theorem 13.3.1]{book} that
we can  find a strongly symmetric
recurrent Borel right
process $Y$
   on $S\cup\{0\}$ with $P^{ x}( T_{0}<\ff)>0$ for all $x\in S$
such that $X$ is the process obtained by killing
$Y$ the first time it hits $0$.  Let   $ L^x_t=\{ L^x_t; t\in R_{+}, x\in S\cup \{0\}\}$  denote the local time
of $Y$. It follows
from the generalized second Ray-Knight Theorem in \cite{fiveauthors}, see
also \cite[Theorem 8.2.2]{book}
  that under $P^{ 0}\times P_{ G}$,
\be
\bigg\{ L^x_{\tau(t)}+{ 1\over 2}{\(G_x\over c_x\)^2};\,x\in S
\bigg\}\stackrel{law}{=}
\bigg\{{ 1\over 2}\({G_x\over c_x}+\sqrt{2t}\)^2;\,x\in S\bigg\}\label{nor.1rqq}
\ee
   for all $t\in R_+$, where $\tau(t)=\inf \{s>0|L_s^0>t\}$,   the inverse local time at zero, and $Y$ and $G$
are independent. Consequently 
\be
\Big\{c_x^2 L^x_{\tau(\al^2/2)}+{ 1\over 2}{G_x ^2};\,x\in S
\Big\}\stackrel{law}{=}
\Big\{{ 1\over 2}\(G_x + c_x \al\)^2;\,x\in S\Big\}\label{nor.1ww}
\ee
for all $\al\in R^1$. (We can extend $\al$ from $R_+$ to $R^1$ because $G$
is symmetric.) $\{c_x^2 L^x_{\tau(\al^2/2)} ;\,x\in S\}$ and $ \{{ 1\over 2}{G_x
^2};\,x\in S\}$ are independent. $G^2$ is infinitely divisible and for all
integers $r\ge 1$
\be 
c_{ \cd}^2 L^{\cdot}_{\tau(\al^2/2)}\stackrel{law}{=}c_{\cd}^2
L^{\cdot}_{\tau(\al^2/(2r)) ,1}+\cdots +c_{ \cd}^2
L^{\cdot}_{\tau(\al^2/(2r)) ,r}\label{1.17}
\ee
where $\{L^{\cdot}_{\tau(\al^2/(2r)) ,j}\}$, $j=1,\ldots,r$ are independent.

  Note that in (\ref{nor.1ww})  we identify the components of the decomposition of $\{\(G_x + c_x \al\)^2;\,x\in S\}$ that mark it as infinitely divisible.

  \medskip	  In Theorem \ref{lem-Bapat} we have necessary and sufficient conditions for   $((G_{1}+c _{1}\al)^{2}, (G_{2}+c _{2}\al)^{2}) $ to be infinitely divisible   for all $ \al\in R^1$. There remains the question, can  $((G_{1}+c _{1}\al)^{2}, (G_{2}+c_{2} \al)^{2}) $  have infinitely divisible squares  for some   $ \al>0$ but not for all $ \al\in R^1$? 
  When we began to investigate this question we hoped that such   points $\al$  do not exist. This would   have finished off the problem of characterizing Gaussian random variables with infinitely divisible squares and, more significantly, by  (\ref{nor.1ww}), would show that when a Gaussian random variable  with non--zero mean has infinitely divisible squares,
it  decomposes into the sum of two independent random variables. The Gaussian random variable itself minus its  mean, and the local time  of   a related Markov process. This would be a very neat result indeed, but it is not true.

For all  Gaussian random variables  $(G_1,G_2)$ in $R^{2}$ and all $c_{1},c_{2}\in R^{1}$ define
 \be
\GG^{2}(c_{1},c_{2},{\al}):=((G_{1}+c_{1} \al)^{2}, (G_{2}+c_{2} \al)^{2}). 
\ee
It follows from Theorem \ref{lem-Bapat}, (for details see \cite[Corollary 1.3, 4.]{IDS}), that $\GG^{2}(c_{1},c_{2},{\al})$ has infinitely divisible squares   for all $ \al\in R^1$ if and only if
\bea 
   \Ga_{1,1}\ge \frac{c_1}{c_2}\Ga_{1,2}\quad\mbox{and}\quad  \Ga_{2,2}\ge \frac{c_2}{c_1}\Ga_{1,2}.\label{1.12q}
 \eea
  If (\ref{1.12q}) does not hold,  we call  $0<\al_0<\ff$    a {\bf critical point} for the infinite divisibility of
  $\GG^{2}(c_{1},c_{2},{\al})$    if    $\GG^{2}(c_{1},c_{2},{\al})$ is infinitely divisible for all $|\al|\le \al_{0}$, and is not infinitely divisible for any $|\al|>\al_{0}$. In this paper we prove the following theorem:

\medskip	
\bt\label{theo-cp}
  For all Gaussian random variables  $(G_1,G_2)$ in $R^{2}$ and all   $(c_{1},c_{2}) \in R^{2}$ for which (\ref{1.12q}) does not hold,  $\GG^{2}(c_{1},c_{2},{\al}) $ has a critical point.   \et

Note that in Theorem \ref{theo-cp} we consider all  $(c_{1},c_{2})\in R^{2}$. It follows from (\ref{1.12q}) that when $EG_{1} G_{2}>0$, then  $\GG^{2}(c_{1},c_{2},{\al})$ has infinitely divisible squares for all $\al\in R^{1}$  only if  $c_{1}c_{2}>0$.  Nevertheless, by Theorem \ref{theo-cp}, even when $c_{1}c_{2}\le 0$,  $ (G_1+{c_{1}\al},G_2+c_{2}\al)$ does  have  infinitely divisible squares for   $|\al|$ sufficiently small.

\medskip	To conclude this Introduction we explain how we approach the problem of showing that $(G+c\al)^2$ is infinitely divisible for all $\al\in R^1$ or only for some $\al\in R ^{1}$. Since we can only prove Theorem \ref{theo-cp} for Gaussian random variables in $R^{2}$ we stick to this case, although a similar analysis applies to $R^{n}$ valued Gaussian random variables. 

 Let $\Ga$ be the covariance matrix of $G$ and 
\be 
\wt\Ga := (I+\Ga\La)^{-1}\Ga=(\Ga^{-1}+\La)^{-1} ,\label{1.12xx}
\ee 
where
\be\La=
     \(\begin{array}{cc}\la_1&0\\
 0&\la_2\end{array}\). 
 \ee
Consider the Laplace transform of  $((G_{1}+c_{1}\al  )^2, (G_2+c_{2}\al)^2)$, 
\bea
\lefteqn{ E_G\(e^{ -( \la_1(G_{1}+c_{1}\al  )^2+\la_2(G_2+c_{2}\al)^2)/2}\)\label{ids.qaspr.2q}}
\\
&& =\frac{1}{(\det \(I+\Ga\La\))^{1/2} }   \exp\(-  {\al ^2\over 2}   \(  c_{1}^2\la_{1} +c_{2}^2\la_{2} -\sum_{i,j=1}^{2}c_{i}c_{j}\la_{i}\la_{j}\wt\Ga_{i,j}\)\) \nn.
\eea
  (See \cite[Lemma 5.2.1] {book}.)
Set $\la_1=t(1-s_1)$ and $\la_2=t(1-s_2)$,  $0\le s_{1},s_{2}\le 1$ and write (\ref{ids.qaspr.2q}) as 
\bea 
   \exp\(U\( s_{1},s_{2},t,\Ga \)+\al^{2}V\( s_{1},s_{2},t, c_{1} ,c_{2},\wt\Ga \)\)\label{1.15}
 \eea
where  $U:=U( s_{1},s_{2},t, \Ga )=-1/ 2\log(\det \(I+\Ga\La\))$. 
Note that \newline $ \exp\(U\( s_{1},s_{2},t,  \Ga \)\)$ is the Laplace transform of $(G_{1}^{2},G_{2}^{2})$, with the change of variables $\la_1=t(1-s_1)$ and $\la_2=t(1-s_2)$. It is easy to see from Theorem \ref{lem-Bapat0} that all two dimensional Gaussian random variables are infinitely divisible. Therefore, for all $t$ sufficiently large, all the coefficients of the  power series expansion of $U$ in $s_{1}$ and $s_{2}$ are positive, except for the constant term. This is a necessary and sufficient condition for a function to be the Laplace transform of an infinitely divisible random variable. See e.g. \cite[Lemma 13.2.2]{book}. 

Now, suppose that  for all $t$ sufficiently large,
$V\( s_{1},s_{2},t ,c_{1} ,c_{2},\wt\Ga \)$ has all the coefficients of  its power series expansion   in $s_{1}$ and $s_{2}$  positive, except for the constant term.
Then the right--hand side of (\ref{ids.qaspr.2q}) is the Laplace transform of two independent infinitely divisible random variables. It is completely obvious that this holds for all $\al\in R^{1}$. 

On the other hand suppose that  for all $t$ sufficiently large, the power series expansion of $V\( s_{1},s_{2},t ,c_{1} ,c_{2},\wt\Ga \)$   in $s_{1}$ and $s_{2}$ has even one  negative coefficient, besides the coefficient of the constant term.  Then for all $\al$ sufficiently large, (\ref{1.15}) is not the Laplace  transform of an infinitely divisible random variable. In other words  $((G_{1}+c_{1}\al  )^2, (G_2+c_{2}\al)^2)$ is not infinitely divisible for all $\al\in R^{1}$. But it may be   infinitely divisible   if $\al$ is small, since the positive coefficients of U may be greater than or equal to $\al^{2}$ times the corresponding negative coefficients of    $ V$. Clearly, if this is true for some $|\al|=\al_{0}>0$, then it is true for all $|\al|\le \al_{0}$.  

The preceding paragraph explains how we prove Theorem \ref{theo-cp}. We consider vectors $((G_{1}+c_{1}\al  )^2, (G_2+c_{2}\al)^2)$ that are not infinitely divisible for all $\al\in R^{1}$, (this is easy to do using (\ref{1.12q})), and show that for $|\al|$ sufficiently small the coefficients in the power series expansion of
\begin{equation}
   U\( s_{1},s_{2},t,\Ga \)+\al^{2}V\( s_{1},s_{2},t,c_{1} ,c_{2},\wt\Ga \)
   \end{equation}
 in $s_{1}$ and $s_{2}$  are positive, except for the constant term. Our proof  only uses   elementary mathematics, although it is quite long and complicated. In the course of the proof we show that the coefficients of the power series expansion of $U$  in $s_{1}$ and $s_{2}$  are positive, except for the constant term. This provides a direct elementary proof of the fact that the Gaussian random variable $(G_{1},G_{2})$ always has infinitely divisible squares.

As we have just stated, and as the reader will see, the  proof of Theorem \ref{theo-cp} is long and complicated.  So far we have not been able to  extend it to apply to Gaussian random variables in $R^{3}$. One hopes for a more sophisticated and much shorter proof of Theorem \ref{theo-cp} that doesn't depend on the dimension of the Gaussian random variable.

 \section{Gaussian squares in $R^2$ and their Laplace transforms }\label{sec-2}

  Let $G=(G_1,G_2)$ be a mean zero Gaussian process
with covariance matrix 
 \be \Ga= 
     \(\begin{array}{cc}a&1\\
 1&b\end{array}\) \label{cov}\ee
where  $ab=d+1>1$, and let  $G+c:=(G_1+c_{1},G_2+c_{2})$, $c_{1},c_{2}\in R^{1}$.   
Note that
\begin{equation}
\det \Ga=d,\label{30.1}
\end{equation}
and
\begin{equation}
\Ga^{-1}={1 \over d} \(\begin{array}{cc}b&-1\\
-1&a\end{array}\)\label{30.1a}.
\end{equation}
Let  
 \be
 \La=
     \(\begin{array}{cc}\la_1&0\\
 0&\la_2\end{array}\). 
 \ee
Then
\begin{equation}
\Ga^{-1}+\La={1 \over d} \(\begin{array}{cc}b+d\la_{1}&-1\\
-1&a+d\la_{2}\end{array}\)\label{30.1axx}
\end{equation}
and
\bea
\wt\Ga&:=&(I+\Ga\La)^{-1}\Ga=(\Ga^{-1}+\La)^{-1}\\
&=&d\,\, \(\begin{array}{cc}b+d\la_{1}&-1\\
-1&a+d\la_{2}\end{array}\)^{-1} \nn\\
&=&\frac{1}{H(a,b,\la_1,\la_2)  } \(\begin{array}{cc}\la_2d+a&1\\
 1&\la_1d+b\end{array}\) \nn
\eea
where 
\bea
H(a,b,\la_1,\la_2)=1+a\la_1+b\la_2+d\la_1\la_2=d\det \(\Ga^{-1}+\La\).\label{2.7}
\eea
 (We use repeatedly the fact that $ab=d+1$.) Note that by (\ref{30.1}) we have that
 \begin{equation}
\det \(I+\Ga\La\)= \det \(\Ga\(\Ga^{-1}+\La\)\)=d\det \(\Ga^{-1}+\La\)=H(a,b,\la_1,\la_2).\label{30.2}
\end{equation}

\begin{lemma} \label{lem-2.1}
 \bea
\lefteqn{E_G\(e^{ -( \la_1(G_{1}+c_{1}  )^2+\la_2(G_2+c_{2}  )^2)/2}\)\label{idsqaspr.2q}}
\\
&& =\frac{1}{(H(a,b,\la_1,\la_2))^{1/2} }   \exp\(- \frac{c_{1} ^2\la_{1}+c_{2}^{2}\la_{2}+\(c_{1}^{2}b+c_{2}^{2}a-2c_{1}c_{2}\) \la_1
 \la_2 }{2H(a,b,\la_1,\la_2)} \) \nn.
\eea
 \end{lemma}

 \Proof  By \cite[Lemma 5.2.1] {book} 
 \bea
 \lefteqn{E_G\(e^{ -( \la_1(G_{1}+c_{1}  )^2+\la_2(G_2+c_{2}  )^2)/2}\)\label{www} }\\
&&\qquad  =\frac{1}{(H(a,b,\la_1,\la_2))^{1/2} } \exp\(-  {1\over 2}\Bigg (
c_{1}^{2}\la_1 + c_{2}^{2}\la_2\right. \nn\\
&&\qquad
 \qquad\left. -\frac{c_{1}^{2}  \la_1^{2} (\la_2d+a)+2c_{1}c_{2} \la_1\la_2+ c_{2}^{2}\la_2^{2}
(\la_1d+b) }{H(a,b,\la_1,\la_2)}\Bigg)\)   .
\nn
\eea
A simple computation shows that 
\begin{eqnarray}
&&\( c_{1}^{2}\la_1 + c_{2}^{2}\la_2 \)H(a,b,\la_1,\la_2)\label{30.3}\\
&&\hspace{1in}-\(c_{1}^{2}  \la_1^{2} (\la_2d+a)+2c_{1}c_{2} \la_1\la_2+ c_{2}^{2}\la_2^{2}
(\la_1d+b)\) 
\nn\\
&&\qquad=c_{1} ^2\la_{1}+c_{2}^{2}\la_{2}+\(c_{1}^{2}b+c_{2}^{2}a-2c_{1}c_{2}\) \la_1
 \la_2 \nn,
\end{eqnarray}
from which we get (\ref{idsqaspr.2q}).\qed

\medskip	

The term   $1/(H(a,b,\la_1,\la_2))^{1/2}$ is
the Laplace transform of 
$(G_1^2,G^2_2)/2$ and  by   \cite[Corollary 1.1, 2.]{IDS}   it is the Laplace
transform of an infinitely divisible random variable.   The exponential term may or may not be a Laplace transform. In
fact, by (\ref{1.12q}), we know it is the Laplace transform of an infinitely divisible
random variable,   for all $\{c_{1}\al,c_{2}\al\}$, for  all $\al\in R^{1}$,  if and only if 
\be
a\ge \frac{c_{1}}{c_{2}}>0\qquad \mbox{and}\qquad b\ge \frac{c_{2}}{c_{1}}>0.\label{4.5}
\ee
  To prove Theorem \ref{theo-cp} we must show that when  (\ref{4.5}) does not hold, there exists    an  $\,0<\al_{0}<\ff$ such that  (\ref{idsqaspr.2q}) is the Laplace transform of an infinitely divisible random variable when $c_{1}$ and $c_{2}$ are replaced by   $c_{1}\al$ and $c_{2}\al$ for any $|\al|\le \al_{0}$.   Actually,  as we see in Section \ref{sec-8}, the general result follows from the consideration of   three cases, 
\begin{enumerate}
\item $c_{1}=c_{2}:=({c,c})$;
\item $c_{1}=-c_{2}:=({c,-c})$ 
\item $(c,0)$. 
\end{enumerate}
This is because if $c_{1}\ne c_{2}$ and neither of them is zero, we can replace $(G_{1},G_{2})$ by   $(G_{1}/|c_{1}|,G_{2}/|c_{2}|)$. Clearly, in this case,  if Theorem \ref{theo-cp} holds for $(G_{1}/|c_{1}|,G_{2}/|c_{2}|)$ it holds for $(G_{1},G_{2})$.  

In these three cases  the  numerator of the fraction in the exponential term on the right--hand side of
 (\ref{idsqaspr.2q}) is
 \begin{enumerate}
\item $ c^{2}\((a+b-2)\la_{1}\la_{2}+\la_{1}+\la_{2}\)$;
\item  $ c^{2}\((a+b+2)\la_{1}\la_{2}+\la_{1}+\la_{2}\)$ 
\item $ c^{2}\(b\la_{1}\la_{2}+\la_{1} \)$. 
\end{enumerate}
Set 
\begin{equation}
 \ga=  a+b-2\qquad\mbox{and}\qquad \rho=a+b+2.
   \end{equation}
 Note unless $\det \Ga=0$, $ab>1$. Since  Theorem \ref{theo-cp} obviously holds
 when $\det \Ga=0$, we can exclude this case from further consideration. Thus we always have $\ga>0$.

 \section{Power series expansion of the   logarithm of the Laplace transform of ${\bf((G_{1}+c)^{2}, (G_{2}+c)^{2})}$   when ${\bf EG_{1} G_{2} =1}$}\label{sec-3}
 
  Bapat's proof of Theorem \ref{lem-Bapat0}   involves  the analysis of a certain power series expansion of the logarithm of the Laplace transform. We need a similar, but more delicate,  analysis. (See Lemma \ref{lem-1.1} below).

Using (\ref{idsqaspr.2q})  and the remarks following Lemma \ref{lem-2.1} we can write  
\bea
\lefteqn{E_G\(e^{ -( \la_1(G_{1}+c )^2+\la_2(G_2+c 
)^2)/2}\)\label{idpr.2}}
\\
&&\qquad  =\exp\(-\frac{1}{2}\log H(a,b,\la_1,\la_2)
\)     \exp\(- \frac{c ^2\(\ga\la_1
 \la_2+\la_1+\la_2\)}{2H(a,b,\la_1,\la_2)} \) 
. 
\nn\\
&&\qquad:  =\exp\( \frac{1}{2}\Big(P(a,b,\la_1,\la_2)+
c ^2 Q(a,b,\la_1,\la_2  ) \Big)\) ,
\nn
\eea
Since $ab=d+1$,   (recall that $d>0$), we have
\begin{eqnarray}
a+b-(d+2)&=&a+{(d+1) \over a}-(d+2)
\label{4.5jss}\\
&=& {1 \over a }\(a^{2}-(d+2)a+(d+1)\)  \nonumber\\
&=& {1 \over a }\(a-(d+1)\)\(a-1\).  \nonumber
\end{eqnarray}
Thus $a+b-(d+2)\le 0$ if and only if $1\le a \le d+1$, which in view of $ab=d+1$ is equivalent to $1\le b\le d+1$.
Consequently (\ref{4.5}) holds  if and only if  $a+b-(d+2)\le 0$. 

Therefore, to show that  $((G_{1}+c)^{2},(G_{2}+c)^{2})$ is infinitely divisible, for some, but not for all,   $c>0$, we must consider $a, b>0$ such that  
\begin{equation}
\ze := a+b-(d+2)>0.\label{4.5k}
\end{equation}
 In the rest of this paper we   assume that (\ref{4.5k}) holds.

\medskip	
Let $\la_1=t(1-s_1)$ and $\la_2=t(1-s_2)$,  $0\le s_{1},s_{2}\le 1$. We consider $P$ and $Q$ as
functions  of and $s_1,s_2 ,t$, and write
\begin{equation}
P(s_1,s_2 ,t):=P(a,b,\la_1,\la_2),\hspace{.3 in}Q(s_1,s_2 ,t ):=Q(a,b,\la_1,\la_2 ).\label{22.1}
\end{equation}

We  expand these in a power series in $s_1,s_2$. 
\begin{equation}
P(s_1,s_2 ,t)=\sum_{j,k=0}^{\ff}P_{j,k}(t)s_{1}^{j}s_{2}^{k},\hspace{.3 in}Q(s_1,s_2 ,t )=\sum_{j,k=0}^{\ff}Q_{j,k}(t)s_{1}^{j}s_{2}^{k}\label{22.2},
\end{equation}
and    set 
\begin{equation}
R(s_1,s_2 ,t,c)=P(s_1,s_2 ,t)+c^{2}Q(s_1,s_2 ,t,c)\label{22.2s}.
\end{equation}
Consequently
\begin{equation}
R(s_1,s_2 ,t,c)=\sum_{j,k=0}^{\ff}R_{j,k}(t,c)s_{1}^{j}s_{2}^{k}\label{22.2t}
\end{equation}
with
\begin{equation}
R_{j,k}(t,c)=P_{j,k}(t)+c^{2}Q_{j,k }(t).\label{22.2u}
\end{equation}

In this section we obtain explicit expressions for $P_{j,k}(t),Q_{j,k}(t)$.
\medskip   
We write
\bea
H(a,b,\la_1,\la_2)&=&{ 1+a\la_1+b\la_2+d\la_1\la_2}\label{3.11}\\
&=&1+at+bt+dt^2 -(at +dt^2 )s_1-(bt+dt^2 )s_2+dt^2 s_1s_2\nn\\
&=&(1+at +bt+dt^2 )(1-\al s_1-\bb s_2+ps_1s_2)\nn\\
&=&(1+at +bt+dt^2 )(1-\al s_1-\bb s_2+\th\al\bb s_1s_2)\nn
\eea
 where
\be
\al={at +dt^2 \over 1+at +bt+dt^2 },\quad \bb={bt+dt^2 \over 1+at +bt+dt^2 },\quad
p={ dt^2 \over 1+at +bt+dt^2 }\label{as1.5q}
\ee
and
 \be
\al\bb=p\(1+d^{-1}+at +bt+dt^2 \over 1+at +bt+dt^2 \):=\frac{p}{\th}\label{as1.5}.
\ee
Note that
\be
1-\th=\frac{d^{-1} }{1+d^{-1}+at +bt+dt^2 }\leq {1 \over d^{2}t^{2}}  .\label{as1.5qq}
\ee

 Using these definitions we have 
\begin{eqnarray}
  P(a,b,\la_1,\la_2 )&=&-\log H(a,b,\la_1,\la_2)
\label{3.10p}\\
& =&-\log (1+at +bt+dt^2 )-\log (1-\al s_1-\bb s_2+\th\al\bb s_1s_2)   \nonumber
\end{eqnarray}
and 
\bea 
&&
 Q(a,b,\la_1,\la_2 )\label{3.10}\\
 &&\qquad=
 -\frac{\ga\la_1 \la_2+\la_1+\la_2 }{ H(a,b,\la_1,\la_2)}\nn\\
 &&\qquad= -{1 \over (1+at +bt+dt^2 )}\,\,\frac{\ga t^2(1-s_1)(1-s_2)+t(2-s_1-s_2) }{ 1-\al s_1-\bb s_2+\th\al\bb s_1s_2}.\nn
\eea

We make some preliminary observations that  enable us to compute the coefficients of the power series expansions of $P(s_1,s_2 ,t)$ and $Q(s_1,s_2 ,t)$. 
  Let $(u_{1},u_{2})\in [0,1)^{2}$, and $\th\in [0,1)$, and assume that 
  \be
u_{1}+u_{2}-\th u_{1}u_{2}<1.\label{2.1}
\ee
(It is clearly greater than zero.)
Let
\be 
{1\over 1-u_1-u_2+\th u_1u_2}:= \sum_{j,k=0}^\ff D_{j,k}u_1^j u_2^k \label{1.1}
\ee 
and
\be 
-\log  (1-u_1-u_2+\th u_1u_2 ):= \sum_{ j,k=0 }^\ff C_{j,k}u_1^j u_2^k .\label{1.2}
\ee 
We give explicit expressions for $C_{j,k}$ and $D_{j,k}$. To begin  we give several equalities that are easy to verify. 

\bl \label{lem-2.1w}   For $u\in [0,1)$
\be
\frac{1}{(1-u)^q}=\sum_{n=0}^\ff{q+n-1\choose n} u^n.
\ee
\el

\Proof To get  ${q+n-1\choose n}$ differentiate $ 1/(1-u)^q$, $n$ times,
divide by $n!$ and set $u=0$. \qed

\medskip We list the following equalities without proof.    
\bl \label{lem-1.2}
\bea
{k\choose p}-{k-1\choose p} &=&{k-1\choose p-1}\label{bineq}\\
{k-1\choose p}&=&{k-p\over k}{k\choose p}\nn\\
{k-1\choose p-1}&=&{p\over k}{k\choose p}\nn\\
{k \choose p}&=&{k-p+1\over p}{k\choose p-1}\nn.
\eea
\el

\bl \label{lem-2.2q} For $0\le j\le k$
\be 
D_{j,k}=\sum_{p=0}^j(1-\th)^p{j\choose p}  {k\choose p} . \label{2.3} 
\ee
Also, $C_{0,0}=0$, $C_{j,0}= 1/j,C_{0,k}= 1/k, j,k\ne 0$, and for $1\le j\le k$
\be 
C_{j,k}=  \sum_{p=0}^{j-1}   (1-\th )^{p}   {j \choose p } {k  \choose p }{(1-\th )\over  p+1}{(j-p)(k-p)\over jk} . \label{2.3h} 
\ee
\el

\Proof Note that  
\be 
{1\over 1-u_1-u_2+\th u_1u_2}=\sum_{n=0}^\ff\(u_1+u_2-\th u_1u_2\)^n .
\ee 
  Writing (\ref{2.1}) in the form $ u_{1}+u_{2}- u_{1}u_{2}+(1-\th) u_{1}u_{2}<1$ we see that it is equivalent to the statement that
   \begin{equation}
   \frac{(1-\th)u_{1 }u_{2}}{(1-u_{1})(1-u_{2})}<1.
   \end{equation}
 We write
\bea
&&{1\over 1-u_1-u_2+\th u_1u_2}\\
&&\qquad = \( (1-u_1) (1-u_2)\(1-{(1-\th )u_1u_2\over (1-u_1) (1-u_2)}\)\)^{-1}
 \nn \\
&&\qquad= {1\over (1-u_1) (1-u_2)}\sum_{p=0}^\ff\({(1-\th ) u_1u_2\over
(1-u_1) (1-u_2)}\)^p\nn\\
&&\qquad= \sum_{p=0}^\ff {(1-\th )^p u_1^p u_2^p\over (1-u_1)^{p+1}
(1-u_2)^{p+1}} \nn.
\eea

Using this series we can determine $D_{j,k}$. Since $j\le k$ it is clear that
each of the first $j+1$ terms in the series immediately above can contribute
a term in $u_1^ju_2^k$. For a given $0\le p\le j$ we get $(1-\th)^p$ times
the   coefficient of  $u_1^{j-p}$ in the power series expansion of
$1/(1-u_1)^{p+1}$ and the  coefficient  of $u_2^{k-p}$ in the power series
expansion of $1/(1-u_2)^{p+1}$. We see from Lemma \ref{lem-2.1w} that they
are ${j\choose j-p}={j\choose  p}$ and ${k\choose k-p}={k\choose p}$
respectively. Thus we get  (\ref{2.3}).

To obtain (\ref{2.3h}) we write  
\bea
&&-\log(1-u_1-u_2+\th u_1u_2)\\
&&\qquad =-\log((1-u_1)(1-u_2)-(1-\th) u_1u_2)\nn\\
&&\qquad=-\log(1-u_1)-\log(1-u_2)-\log\(1-{(1-\th) u_1u_2\over (1-u_1)(1-u_2)}\)\nn\\
&&\qquad =\sum_{n=1}^\ff \frac{u_1^n}{n}+\sum_{n=1}^\ff \frac{u_2^n}{n}+\sum_{p=1}^\ff 
\frac{1}{p}\({(1-\th ) u_1u_2\over
(1-u_1) (1-u_2)}\)^p\nn.
\eea
This gives us $C_{j,0}$ and $C_{0,k}$ and, similar  to the computation of $D_{j,k}$,  we can use the last series above to see that   
 \bea
C_{j,k}&=& \sum_{p=1}^j {(1-\th )^p \over
p } {j-1\choose p-1} {k-1\choose p-1} \\
&=& \sum_{p=0}^{j-1} { (1-\th )^{p+1}   \over
p+1 } {j-1\choose p } {k-1 \choose p }\nn\\
&=& \sum_{p=0}^{j-1} { (1-\th )^{p+1}\over
p+1 } {j-1\choose p } {k-1 \choose p }\nn\\
 &=& \sum_{p=0}^{j-1}   (1-\th )^{p}   {j \choose p } {k  \choose p }{(1-\th )\over  p+1}{(j-p)(k-p)\over jk} \nn.
\eea
This gives us (\ref{2.3h}). \qed

  Set  \begin{equation}
 \wt t={1 \over \sqrt{1-\th}}.\label{2.51a}
 \end{equation}
 By (\ref{as1.5qq}) we have that 
  \begin{equation}  
dt\leq  \wt t \leq dt+  2,\label{2.51a1}
 \end{equation}
  for all $t$ sufficiently large.

\bl\label{lem-P}  For all $t$ sufficiently large,  $P_{j,0}(t)= \al^{j}/{j}$,   $P_{0,k}(t)= \bb^{k}/{k}$ and for  all $1\le j\le k$  
\be 
P_{j,k}(t)={  \al^{j} \bb^{k} \over \wt t^{2}}\sum_{p=0}^{j-1}    \wt t^{-2p}   {j \choose p } {k  \choose p }{1\over  p+1}{(j-p)(k-p)\over jk} . \label{2.3hb} 
\ee
 \el
(See (\ref{22.2}).)

\medskip	
\Proof Note that
since $0<\al,\bb,\th<1$, 
\bea
\al s_1+\bb s_2-\th\al\bb s_1s_2&\le& \al  +\bb  -\th\al\bb\label{4.5l}\\
&=& \al  +\bb  -\al\bb+(1-\th)\al\bb\nn\\
&=& -(1-\al)(1-\bb) +\frac{\al\bb}{d^{2}t^{2}}+1+O(1/t^{3})\nn\\
&=& -\frac{ab}{d^{2}t^{2}} +\frac{1}{d^{2}t^{2}}+1+O(1/t^{3})\nn\\
&=&1-\frac{1}{d^{2}t^{2}} +O(1/t^{3})\nn.
\eea
 Consequently, for all $t$ sufficiently large 
 \begin{equation}
 0\le \al s_1+\bb s_2-\th\al\bb s_1s_2<1.\label{4.5m}
 \end{equation}
 Therefore, by  (\ref{1.2})\be 
-\log  (1-\al u_1-\bb u_2+\al\bb\th u_1u_2 )= \sum_{ j,k=0 }^\ff   \al^{j} \bb^{k} C_{j,k}u_1^j u_2^k .\label{1.2x}
\ee 
The  lemma now  follows from Lemma \ref{lem-2.2q} and (\ref{2.51a}).
\qed

 Set  \begin{equation}
 \bar t= \sqrt{1+at+bt+d t^2}\label{2.51asr}
 \end{equation}
and note that  
  \begin{equation}  
  d^{1/2} t\leq  \bar t \leq d^{1/2}t+2\label{2.51asr1}
 \end{equation}
  for all $t$ sufficiently large.

\bl\label{lem-2.1j}  For  all $t$ sufficiently large, and $j,k\ge 1$  
 \be 
Q_{j,0}(t)= -\frac{ ( \ga t^2(\al-1) +t(2\al-1)
 )}{{\bar {t}}^2} \al^{j-1} \label{3.29o}
 \ee 
and
 \be 
Q_{0,k}(t)= -\frac{ ( \ga t^2(\bb-1) +t(2\bb-1)
 )}{\bar t^2} \bb^{j-1} \nn.\label{3.30o}
 \ee 
Furthermore,   for all $t$ sufficiently large and for all   $1\le j\le k$
\bea
\lefteqn{  Q_{j,k}(t)\label{2.22js}}\\&&={\al^{j-1} \bb^{k-1}\over \bar t^2}\sum_{p=0}^{j}\wt t^{-2p}{j\choose p}{k\choose p}\nn\\
&&\quad\(-\ga t^2\((1-\al)(1-\bb)-{(1-\bb)p\over j}-{(1-\al)p\over k}+\wt t^{-2}{j-p\over j}{k-p\over k}
\)\right.
\nn\\
&&\qquad\left.+t\(\al\(1-\bb-\frac{p}{k}\)+\bb\(1-\al-{p\over j}\) \)\)\nn.
    \eea  
\el

\Proof
 It follows from (\ref{4.5m}) that for $0\le s_1,s_2\le 1$   
\be
\frac{1 }{ 1+a\la_1+b\la_2+d\la_1\la_2}=\frac{1}{\bar t^2}\sum_{n=0}^\ff
(\al s_1+\bb s_2-\th\al\bb s_1s_2)^n.\label{as1.1}
\ee
 Using this along with (\ref{3.10})   we see that 
\bea
Q(s_1,s_2,t )
&=&-\frac{ (\ga t^2(1-s_1)(1-s_2)+t(2-s_1-s_2))}{\bar t^2}\label{as1.3}\\
&&\qquad
\cdot \sum_{n=0}^\ff  (\al s_1+\bb s_2-\th\al\bb s_1s_2)^n\nn\\
&=&-\frac{ (\ga t^2+2t-(\ga t^2+t)
s_1-(\ga t^2+t)s_2   +\ga t^2 s_1s_2)}{\bar t^2}\nn\\ &&\qquad
\cdot\sum_{n=0}^\ff  (\al s_1+\bb s_2-\th\al\bb s_1s_2)^n\nn .
\eea

  Consequently
 \be 
Q_{j,0}(t )=-\frac{ ((\ga t^2+2t)\al-(\ga t^2+t)
 )}{\bar t^2} \al^{j-1},
 \ee 
from which we get (\ref{3.29o}), and 
 \be 
Q_{0,k}(t )=-\frac{ ((\ga t^2+2t)\bb-(\ga t^2+t)
 )}{\bar t^2} \bb^{j-1},
 \ee 
 from which we get (\ref{3.30o}).

  To obtain (\ref{2.22js}) we use (\ref{as1.3}), and the terms   of $D_{j,k}$ defined in (\ref{1.1}),  to see that   
\bea
&&\QQ_{j,k}:= {\bar t^2\,Q_{j,k}\over  \al^{j-1} \bb^{k-1}}\label{1.20}\\
& &\qquad\, \,=-\ga t^2\( D_{j,k}\al\bb -D_{j,k-1}\al-D_{j-1,k}\bb+ D_{j-1,k-1}\)\nn\\
& &\qquad \qquad-t\( 2D_{j,k}\al\bb  -D_{j,k-1}\al-D_{j-1,k}\bb \)\nn.
\eea
  Using (\ref{1.20}),  (\ref{2.3}) and Lemma \ref{lem-1.2} we see that 
\bea
\QQ_{j,k}&=&\sum_{p=0}^{j}\wt t^{-2p}{j\choose p}{k\choose p}\label{2.22}\\
&&\quad\(-\ga t^2\(\al\bb-\al{k-p\over k}-\bb{j-p\over j}+{j-p\over j}{k-p\over k}
\)\right.
\nn\\
&&\qquad\left.-t\(2\al\bb-\al{k-p\over k}-\bb{j-p\over j} 
\)\)\nn .
    \eea
 Consequently,  for all $t$ sufficiently large, for all $1\le j\le k$  
\bea
Q_{j,k}(t )&=&{ \al^{j -1} \bb^{k-1}\over \bar t^2}\sum_{p=0}^{j}\wt t^{-2p}{j\choose p}{k\choose p}\label{2.22j}\\
&&\quad\(-\ga t^2\(\al\bb-\al{k-p\over k}-\bb{j-p\over j}+{j-p\over j}{k-p\over k}
\)\right.
\nn\\
&&\qquad\left.-t\(2\al\bb-\al{k-p\over k}-\bb{j-p\over j} 
\)\)\nn .
    \eea  
 
  Consider  (\ref{2.22j}).
We write
\bea
&&\al\bb-\al{k-p\over k}-\bb{j-p\over j}+{j-p\over j}{k-p\over k}\label{2.23s}\\
&&\qquad=(1-\al)(1-\bb)-{(1-\bb)p\over j}-{(1-\al)p\over k}+{p^2\over jk},\nn
\eea
\be
2\al\bb-\al{k-p\over k}-\bb{j-p\over j}=-\al\(1-\bb-\frac{p}{k}\)-\bb\(1-\al-{p\over j}\),\label{2.24s}
\ee
to obtain  
\bea
Q_{j,k}(t )\label{2.22jt}\\&=&{  \al^{j-1} \bb^{k-1}\over \bar t^2}\sum_{p=0}^{j}\wt t^{-2p}{j\choose p}{k\choose p}\nn\\
&&\quad\(-\ga t^2\((1-\al)(1-\bb)-{(1-\bb)p\over j}-{(1-\al)p\over k}+{p^2\over jk}
\)\right.
\nn\\
&&\qquad\left.+t\(\al\(1-\bb-\frac{p}{k}\)+\bb\(1-\al-{p\over j}\) \)\)\nn.
    \eea 
   Note that for  $1\le q\le j$
 \bea
 &&\wt t^{-2q}{j\choose q}{k\choose q}\({ q^2\over jk}\)\\
 &&\qquad=\wt t^{-2(q-1)}{j\choose q-1}{k\choose q-1}  {j-(q-1)\over j}{k-(q-1)\over k}\wt t^{-2}.\nn
 \eea
 Therefore, for each $1\le p\le j$ we incorporate the term in $ p^2/ jk$
in (\ref{2.22jt}) into the  preceding term in the series in (\ref{2.22jt}) to get 
 (\ref{2.22js}).    (Note that we can not add anything to the $p=j$ term. The expression in (\ref{2.22js}) reflects this fact since $ {j-p\over j}{k-p\over k}=0$ when $p=j$.)
\qed

\section{ A sufficient condition for a vector in $R^{2}$ to be infinitely divisible.}\label{sec-4}

We present a sufficient  condition for   a random vector in $R^{2}$ to be to be infinitely divisible, and show how it simplifies the task of showing that $((G_{1}+c)^{2},(G_{2}+c)^{2})$    is infinitely divisible.

\begin{lemma} \label{lem-1.1}  Let $\psi : (R_{+})^{2}\to  (R_{+})^{2}$ be a continuous function  with $\psi (0,0)=1$. Let   ${\bf s}\in  [0,1]^{2}$ and suppose that for all $t>0$ sufficiently large, $\log \psi(t(1-s_{1}), t(1-s_{2}))$ has a power series expansion at $\bf{s=0}$ given by  
 \begin{equation}
 \phi (t;s_{1},s_{2})=  \sum_{j,k=0}^{\ff} b_{j,k}(t)s_{1}^{j}s^{k}_{2}.\label{1.2sr}
   \end{equation}

  Suppose also that there exist an increasing sequence of  finite subsets   $\NN_{i}\subseteq   \mathbf{N^{2}}$, $i\ge 1$, with $\bigcup_{i=1}^{\ff}\NN_{i}= \mathbf{N^{2}}$, and  a sequence $ t_{i}\to\ff $, $i\ge 1$, such that $b_{j,k}(t_{i})\ge 0$  for all   $(j,k)\in \NN_{i}/(0,0)$ and  
  \begin{equation}
   \lim_{ i\to \ff}\sum_{j,k\notin N_{i}\cup (0,0)}^{\ff} |b_{j,k}(t_{i})|=0.\label{1.2s}
   \end{equation}
 Then  $\psi(\la_{1},\la_{2})$ is the  Laplace transform of an infinitely divisible random variable on $(R_{+})^{2}$.
 \end{lemma}

 \Proof    It is clear from (\ref{1.2s}) that the power series in (\ref{1.2sr}) converges absolutely for all ${\bf s}\in  [0,1]^{2}$.
 Let
    \begin{equation}
 \phi_{i}(t_{i};s_{1},s_{2})=   b_{0,0}(t_{i}) + \sum_{j,k\in N_{i}/(0,0) } b_{j,k}(t_{i})s_{1}^{j}s^{k}_{2}.
   \end{equation} 
   Set  
   \begin{equation}
   \Psi_{i}\(t_{i}; e^{-\la_{1}/t_{i}},  e^{-\la_{2}/t_{i}}\)=\exp\(  \phi_{i} \(t_{i}; e^{-\la_{1}/t_{i}},  e^{-\la_{2}/t_{i}}\) \).
   \end{equation}
   We show that for each  $(\la_{1},\la_{2})\in (R_{+})^{2}$ 
       \begin{equation}
   \lim_{i\to\ff} \Psi_{i}\(t_{i}; e^{-\la_{1}/t_{i}},  e^{-\la_{2}/t_{i}}\)=\psi (\la_{1},\la_{2}).\label{1.6}
   \end{equation}
   
   As we point out in  \cite[page 565]{book}, $    \Psi_{i}\(t_{i}; e^{-\la_{1}/t_{i}},  e^{-\la_{2}/t_{i}}\)$ is the Laplace transform of a discrete    measure. It then follows from the continuity theorem and the fact that $\psi (0,0)=1$ that
$\psi(\la_{1},\la_{2})$ 
  is the Laplace transform   of a random variable. 
 Furthermore repeating this argument with $ \phi_{i} (t_{i};s_{1},s_{2})$ replaced by $  \phi_{i} (t_{i};s_{1},s_{2})/n$ shows that $\psi^{1/n}(\la_{1},\la_{2}) $ is the Laplace transform of a random variable. This shows that $\psi(\la_{1},\la_{2})$ is the  Laplace transform of an infinitely divisible random variable on $(R_{+})^{2}$.
  
Let
 \begin{equation}
\de_{i}:=\Big |\psi(t_{i}(1-e^{-\la_{1}/t_{i}}), t_{i}(1-e^{-\la_{2}/t_{i}}))-\psi(\la_{1},\la_{2})\Big|.\label{1.7}
   \end{equation} 
   Clearly $    \lim_{i\to\ff}\de_{i}=0$.  
By (\ref{1.2s})  
\bea
\lefteqn{\Big| \psi (t_{i}(1-e^{-\la_{1}/t_{i}}), t_{i}(1-e^{-\la_{2}/t_{i}})) - \exp \( \phi_{i}(t_{i};e^{-\la_{1}/t_{i}},e^{-\la_{2}/t_{i}})\)\Big|} \nn \\
&&=\Big| \exp \( b_{0,0}(t_{i})+ \sum_{(j,k)\neq(0,0) } b_{j,k}(t_{i}) e^{-j\la_{1}/t_{i}}e^{-k\la_{2}/t_{i}}\)\\
&&\hspace{.1 in}- \exp  \( b_{0,0}(t_{i})+ \sum_{j,k\in N_{i}/(0,0) } b_{j,k}(t_{i}) e^{-j\la_{1}/t_{i}}e^{-k\la_{2}/t_{i}}\)\Big|\nn\\
&&= \psi(t_{i}(1-e^{-\la_{1}/t_{i}}), t_{i}(1-e^{-\la_{2}/t_{i}}))\nn\\
&&\hspace{1 in}  \Big|1-\exp \Big( -\sum_{j,k\notin  \NN_i\cup (0,0)}^{\ff} b_{j,k}(t_{i}) e^{-j\la_{1}/t_{i}}e^{-k\la_{2}/t_{i}}\Big)\Big| \nn\\
&&= \ep_{i} \psi (t_{i}(1-e^{-\la_{1}/t_{i}}), t_{i}(1-e^{-\la_{2}/t_{i}})) ,\nn
\eea
where
\begin{eqnarray}
\ep_{i} &:=& \Big|\Big(1-\exp \Big( -\sum_{j,k\notin N_{i}\cup (0,0)}^{\ff} b_{j,k}(t_{i}) e^{-j\la_{1}/t_{i}}e^{-k\la_{2}/t_{i}}\Big)\Big)\Big|.
\label{1.10j}
\end{eqnarray}
Note that by (\ref{1.2s})
\begin{equation}
\lim_{i\rar\ff}\ep_{i}=0.\label{1.10k}
\end{equation}
 Therefore, by the triangle inequality,  (\ref{1.7})  and (\ref{1.10k})
     \bea
&&\Big| \exp \( \phi_{i} (t_{i};e^{-\la_{1}/t_{i}},e^{-\la_{2}/t_{i}})\) -  \psi(\la_1,\la_{2})\Big|\label{1.12}\\
&&\qquad\leq \,\ep_{i} \psi (t_{i}(1-e^{-\la_{1}/t_{i}}), t_{i}(1-e^{-\la_{2}/t_{i}}))+  \de_{i} .\nn
\eea   
 Using     (\ref{1.7}) we see that this is 
 \bea
  &&   \le  \ep_{i}\(  \psi(\la_{1},\la_{2}) +\de_{i}\)+\de_{i}.  \nn \eea
  Thus we justify (\ref{1.6}) and the  paragraph following it. \qed

  \begin{remark} {\rm In \cite[Lemma 13.2.2]{book} we present the well known result that the conclusion of Lemma  \ref{lem-1.1} holds when $\log \psi(t(1-s_{1}), t(1-s_{2}))$ has a power series expansion at $\bf{s=0}$ with all its coefficients, except for the coefficient of the constant term, are positive. Lemma  \ref{lem-1.1} is useful because it allows us to only verify this condition for a subset of these coefficients,   (depending on $t$). 
 }\end{remark}
 
 The following lemma enables us to apply Lemma \ref{lem-1.1}. 
     
    \bl\label{lem-cutoff}
For any  $c_{3 }>0$  there exists a constant   $B=B(\ga, d,c_{3})$  for which
    \begin{equation}
    \mathcal{N}_{t}=\{(j,k)\,|\,\sqrt{jk}\leq  Bt\log t\}\label{2.22f}
    \end{equation}
  has the property that
      \begin{equation}
    \lim_{t\to\ff}\sum_{(j,k)\notin   \mathcal{N}_{t}}|R_{j,k}(t,c)|=0\label{2.22kr},
    \end{equation} 
    uniformly in $|c|\leq c_{3}$.
    \el

    \begin{remark}\label{rem-4.2} {\rm It follows from Lemmas \ref{lem-1.1} and \ref{lem-cutoff} that in order to prove Theorem \ref{theo-cp} we need that only show   that we can find   a $c_{0}>0$, such that    \begin{equation}
   R_{j,k}(t,c_{0})\ge 0\qquad \mbox{for all }\quad \sqrt{jk}\leq Bt\log t,
   \end{equation}
for any constant $B$,  for all  $t$ sufficiently large,   
(except for $R_{0,0}(t)$).
 }\end{remark}
 
Before proving   Lemma \ref{lem-cutoff} we   establish the following bounds

\begin{lemma} \label{2.4}
\be
{j\choose p}\le \({ej\over p}\)^{p}\label{4.14}
\ee
and
\be 
\frac1{\tau^{2p}}{j\choose p}  {k\choose p}\le  \({e\sqrt{jk}\over p\, \tau}\)^{2p}\le \exp\({2\sqrt{jk}\over \tau}\) \label{2.12}
 \ee

 \end{lemma}
 
 \Proof
 It is clear that 
 \begin{equation}
   {j\choose p}\le   {j^{p}\over p!}. 
   \end{equation}
   Therefore to prove (\ref{2.4}) we need only show that
   \begin{equation}
 p!\({e\over p}\)^{p}\ge 1.\label{2.13}
   \end{equation}
  In \cite[page 42 ]{Feller1}, Feller shows that $ p!\({e/ p}\)^{p}$ is increasing in $p$. Since it is equal to $e$ when $p=1$, (and 1 when $p=0$), we get (\ref{2.13}).
  
  The first inequality in (\ref{2.12}) follows from (\ref{4.14}) the next one is obtained by maximizing the middle term with respect to $p$.
  \qed

 \medskip	
 
 \noindent   
{\bf  Proof of Lemma \ref{lem-cutoff}} By (\ref{22.2u}) and Lemmas \ref{lem-P} and \ref{lem-2.1j}
we see that for all $t$ sufficiently large, for all $1\le j\le k$, 
    \begin{equation}
   |R_{j,k}(t,c )|\le C\al^{j}\bb^{k} \sum_{p=0}^{j}\wt t^{-2p}{j\choose p}{k\choose p}
      \end{equation}
where $C$ depends on $c,\ga$ and $d$ but not on $j$, $k$,  $t$ or   $\wt t$. Furthermore, $C$ is bounded for all $|c|\le T$, for any finite number $T$.
(We also use the fact that $\lim_{t\to\ff}\al\bb=1$.)

  For any $\de>0$,   for  $t$ sufficiently large,  
 \begin{equation}
\al={at +dt^2 \over 1+at +bt+dt^2 }= 1- {1+bt \over 1+at +bt+dt^2 }\leq 1-{(1-\de)b \over dt}\leq e^{-(1-\de)b/(dt)} \label{2.22g}
 \end{equation}
 and  
  \begin{equation}
\bb={bt+dt^2 \over 1+at +bt+dt^2 }= 1- {1+at \over 1+at +bt+dt^2 }\leq 1-{(1-\de)a \over d t}\leq e^{-(1-\de)a/(dt)}.\label{2.22h}
 \end{equation}

   Using these estimates along with  (\ref{2.51a1}) we see that  for  all $t$ sufficiently large, for all $1\le j\le k$, 
 \begin{equation}
|R_{j,k}(t,c)|\leq  Ce^{-k(1-\de)a/(dt)}e^{-j(1-\de)b/(dt)}\,\sum_{p=0}^{j }{1 \over(d t)^{2p}}{j\choose p}{k\choose p}\label{2.22i}
 \end{equation}
  uniformly in          $|c|\leq c_{3}$.  
 
 Suppose that $\sqrt{jk}/(dt)=n$. Then 
   \bea
   &&
e^{-k(1-\de)a/(dt)}e^{-j(1-\de)b/(dt)}\label{2.22j1}\\
&&\qquad=\exp\(-(1-\de)\( a\sqrt{k/j}+ b\sqrt{j/k}\)n\)\nn\\
   & &\qquad\le\exp\(-2 (1-\de) \sqrt{ab} n\)\nn\\
    & &\qquad=\exp\(- 2(1-\de)\sqrt{d+1}n\)\nn,
   \eea
   where, for the inequality we take the minimum of $a\th+b/\th$ and for the equality we use the fact that $ab=d+1$. 
   Combined with (\ref{2.12}) this shows that when $\sqrt{jk}/(dt)=n$
   \bea
|R_{j,k}(t,c)| 
   &\le& C j  \exp\(- 2((1-\de)\sqrt{d+1}-1)n\).\label{2.22jx}
   \eea    
Let $A_{n}=\{(j,k)| n\le \sqrt{jk}/(dt)<n+1\}$. Then   for any $M$
   \begin{equation}
     \sum_{  \sqrt{jk}/(dt) \geq M}|R_{j,k}(t,c)|=  \sum_{n=M}^{\ff} \sum_{(j,k)\in A_{n}} |R_{j,k}(t,c)|\label{2.22k} .
   \end{equation}
We note that the cardinality of $A_{n}$ is less than  
\be 3nd^{2}t^{2}\log ( (n+1)d t) \label{dfd}.\ee
  This is because $(j,k)\in A_{n}$ implies that 
  \begin{equation}
   {n^{2}d^{2}t^{2}\over j}\le k\le  {(n+1)^{2}d^{2}t^{2}\over j}\quad\mbox{and}\quad j\le (n+1)dt 
   \end{equation}
  and summing on $ j$  we get (\ref{dfd}).

 It follows from (\ref{2.22j}) and (\ref{2.22k}) that 
  \bea
    && \sum_{  \sqrt{jk}/(dt) \geq M }|R_{j,k}(t,c)| \label{2.22jx} \\
   & &\qquad\le C ( dt)^{4} \sum_{n=M}^{\ff}n^{4}
    \exp\(- 2((1-\de)\sqrt{d+1}-1)n\)\nn  \\
    & &\qquad\le C {( dt)^{4}\over (1-\de)\sqrt{d+1}-1)}  M^{4} \exp\(- 2((1-\de)\sqrt{d+1}-1)M\)\nn.
   \eea
Clearly, there exists a constant $B$ such that when  $M=B\log  t$,  this last term is  $o(1)$ as $t\to\ff$. \qed

   \section{ Proof of Theorem  \ref{theo-cp}  when {$\bf (c_{1},c_{2})=(c,c)$} and ${\bf EG_{1} G_{2} >0}$}\label{sec-5}

  In this section we prove Theorem \ref{theo-cp}  in case 1. and for $  EG_{1} G_{2} >0 $,   by   establishing the positivity conditions on the coefficients $R_{j,k}(t,c) $, (when $  EG_{1} G_{2} =1$),  as  discussed in Remark \ref{rem-4.2}. We pass to the case $  EG_{1} G_{2}>0$ on page \pageref{sec5page}.

To proceed we need several estimates of parameters we are dealing with as $t\to\ff$. They follow   from the definitions in (\ref{as1.5q})--(\ref{as1.5qq}).
 
 \bl \label{lem-2.2}As $t\to\ff$
\bea
1-\al&=&{b\over dt}- {1+b^{2}\over (dt)^2 }+O(t^{-3})\label{3.18}\\
1-\bb&=&{a\over dt}- {1+a^{2}\over (dt)^2 }+O(t^{-3})\nn\\
(1-\al)(1-\bb)&=&{
d+1\over (dt)^2}- {a(1+b^{2})+b(1+a^{2})\over (dt)^3 }+O(t^{-4})\nn\\
\wt t^{-2}=1-\th &=&{1\over (dt)^2}-{a+b\over (dt)^3}+O(t^{-4})\nn\\
\al^{j}&=&e^{-bj/(dt)+O(j^{2}/t^{2})}\nn\\
\bb^{k}&=&e^{-ak/(dt)+O(k^{2}/t^{2})}\nn.
\eea

Also
\bea
-(d+2)\ga+d(a+b)&=&2((d+2)-(a+b))=-2\ze\label{2.28}\\
a\ga-d&=&(a-1)^2\nn\\
b\ga-d&=&(b-1)^2\nn.
\eea
\el
 
 \Proof 
 \begin{eqnarray}
 1-\al&=&{1 +bt \over 1+at +bt+dt^2 }
 \label{3.18a}\\
 &=&  {1 +bt \over dt^2 }\,\,{1  \over 1+a(dt)^{-1} +b(dt)^{-1}+d^{-1}t^{-2} } \nonumber\\
 &=&  {1 +bt \over dt^2 }\,(1-a(dt)^{-1} -b(dt)^{-1}+O(t^{-2}))\nonumber\\
 &=& {b\over dt}+ {d -b(a+b)\over d^{2}t^2 }+O(t^{-3})\nonumber\\
 &=& {b\over dt}- {1+b^{2}\over d^{2}t^2 }+O(t^{-3})\nonumber
 \end{eqnarray}
 The rest of the lemma follows similarly.
 \qed

  \medskip	  \noindent  
{\bf Proof of Theorem \ref{theo-cp}      when    ${\bf c_{1}=c_{2}=c}  $}.     To begin let note that it is easy to see from Lemma \ref{lem-P},  that $P_{j,k}(t)\ge 0$ for all $0\le j,k<\ff$, with the exception of $P_{0,0}(t)$. This must be the case because $\exp(P(a,b,\la_{1},\la_{2}))$ is the Laplace transform of an infinitely divisible random variable, as we remark following the proof of Lemma \ref{lem-2.1}.

By (\ref{3.29o})  
\bea
Q_{j,0}(t)  &=&-\frac{ ( \ga t^2(\al-1) +t(2\al-1)
 )}{\bar t^2} \al^{j-1} \label{5.4}
 \\
 &=&-\frac{ (( -\ga b+1)t +\ga(1+b^{2})-2b+O(1/t) )
 }{\bar t^2} \al^{j-1} \nn\\
 &=& { \((b-1)^{2}t+2b-\ga(1+b^{2}) +O(1/t)\)\over {\bar t^2}}\al^{j-1}\nn\\
  &=&\( {  (b-1)^{2} \over dt}+O(1/t^{2})\)\al^{j-1}\nn.
 \eea 

Similarly
 \bea
Q_{0,k}(t )  &=&  \( {  (a-1)^{2} \over d t}+O(1/t^{2})\)\al^{k-1} .\label{3.23}
 \eea 
 Thus  we see that   there exists a $t_{1}$ sufficiently large such that for all $t\geq t_{1}$, $R_{j,0}(t,c)$ and $R_{0,k}(t,c)$ are both positive for all $j,k\ge 1$.

\medskip	
We now examine $R_{j,k}(t,c) $ for $j\wedge k\ge 1$, $j\le k$.
We write
   \begin{equation}
   R_{j,k}(t,c)=\sum_{p=0}^{j}   R_{j,k,p}(t,c)\label{2.4c}.
    \end{equation}
 Using (\ref{2.22js}) and (\ref{2.3hb}) we see that
    \begin{eqnarray}
    \lefteqn{{R_{j,k,p}(t,c) \over  \al^{j-1}\bb^{k-1}}    \label{2.4f} }\\
    &&= {\al\bb \wt t^{-2p} \over \wt t^{2}}  {j \choose p } {k  \choose p }{1\over  p+1}{(j-p)(k-p)\over jk}
  +{c^2\wt t^{-2p}\over \bar t^2}{j\choose p}{k\choose p}\nn\\
&&\quad\(-\ga t^2\((1-\al)(1-\bb)-{(1-\bb)p\over j}-{(1-\al)p\over k}+\wt t^{-2}{j-p\over j}{k-p\over k}
\)\right. \nonumber\\
&&\quad\qquad\left.+t\(\al\(1-\bb-\frac{p}{k}\)+\bb\(1-\al-{p\over j}\) \)\).\nn
    \end{eqnarray}
When $p=0$
we get
      \begin{eqnarray}
    &&{R_{j,k,0}(t,c) \over  \al^{j-1}\bb^{k-1}}={\al\bb \over \wt t^{2}}  
    +{c^2\over \bar t^2} \label{2.4f0}\\
&&\quad\quad\(-\ga t^2\((1-\al )(1-\bb )+\wt t^{-2}\)+t\(\al(1-\bb )+\bb(1-\al ) \) \)\nn
    \end{eqnarray}
   which is independent of $j,k$.   Using Lemma \ref{lem-2.2} we see that  
   \begin{eqnarray}
   &&-\ga t^2\((1-\al )(1-\bb )+\wt t^{-2}\)+t\(\al(1-\bb )+\bb(1-\al ) \)
   \label{2.4f1}\\ 
   &&\qquad=-{(d+2) \over d^{2}}\ga +{a+b \over d}+O\( {1 \over t} \) \nonumber\\
   &&\qquad={-(d+2)\ga+d(a+b)  \over d^{2}}+O\( {1 \over t} \) \nonumber\\
   &&\qquad  ={-2\ze\over d^{2}}+O\( {1 \over t} \) \nonumber.
   \end{eqnarray}
    Using  this and Lemma \ref{lem-2.2} again we get 
   \begin{equation}
  {R_{j,k,0}(t,c) \over  \al^{j-1}\bb^{k-1}}= {1-2c ^2(\ze/d)+O\( 1/t \)\over d^{2} t^2}\label{2.4fx}
   \end{equation}
   where the $O\( 1/t\)$ term is independent of $j$ and $k$.

   \medskip	  
   
   We now simplify the expression of the other coefficients   $R_{j,k,p}(t,c)$, $1\le p\le j$.    Set  
    \begin{equation}
{R_{j,k,p}(t,c)\over \al^{j-1}\bb^{k-1}}=\wt t^{-2p}{j\choose p}{k\choose p}\({1 \over \wt t^2}F_{j,k,p}(t)+{c^2\over \bar t^2} A_{j,k,p}(t)\)\label{5.10}
   \end{equation}
   where
 \begin{equation}
  F_{j,k,p}(t)  =\al\bb    {1\over  p+1}{(j-p)(k-p)\over jk}\label{5.11}
   \end{equation} 
   and 
    \begin{eqnarray}
     \lefteqn{A_{j,k,p}(t)\nn}\\
   &&=-\ga t^2  \((1-\al)(1-\bb)-{(1-\bb)p\over j}-{(1-\al)p\over k}+\wt t^{-2}{j-p\over j}{k-p\over k}\)
 \nn\\
   && \qquad + t    \(\al\(1-\bb-\frac{p}{k}\)+\bb\(1-\al-{p\over j}\) \)   \label{5.12}.
   \eea 
 
     Using Lemma \ref{lem-2.2}  we have  
   \begin{eqnarray}
   &&-\ga t^2  \((1-\al)(1-\bb)-{(1-\bb)p\over j}-{(1-\al)p\over k}+\wt t^{-2}{j-p\over j}{k-p\over k}\)
 \nonumber\\
 &&=-{\ga \over d^{2}}\(d+1  -{(bdt-(1+b^{2}))p\over k}-{(adt-(1+a^{2}))p\over j} +{j-p\over j}{k-p\over k}\)\nn\\
   &&\hspace{4 in}+O\( {1 \over t} \)
   \nn\\
   && ={\ga \over d^{2}}\(-(d+2)  +{(bdt-b^{2})p\over k}+{(adt-a^{2})p\over j} -{p^{2}\over jk}\)+O\( {1 \over t} \) \label{2.4f2}\\
   && ={\ga \over d^{2}}\(-(d+2)   +{b(dt-b)p\over k}+{a(dt-a)p\over j} -{p^{2}\over jk}\)+O\( {1 \over t} \) \nn
   \end{eqnarray}
   and
   \begin{eqnarray}
   &&  t    \(\al\(1-\bb-\frac{p}{k}\)+\bb\(1-\al-{p\over j}\) \)
   \label{2.4f3}\\
   &&\qquad= {a \over d}\(1+{p \over j}\)+ {b \over d}\(1+{p \over k}\)-{  pt\over j}-{ pt\over k} +O\( {1 \over t} \)  \nonumber\\
   &&\qquad= {1 \over d} \( a+b-{  p(dt-a)\over j}-{ p(dt-b)\over k}\) +O\( {1 \over t} \)  \nonumber.
   \end{eqnarray}
   
   In (\ref{2.4f2}) and (\ref{2.4f3}) the expressions $O\( {1 / t} \)$
are not necessarily the same from line to line. Nevertheless, it is important to note that they are   independent of $p$, $j$ and $k$. That is there exists an $M>0$ such that all terms given as   $O\( {1 /t} \)$ in (\ref{2.4f2}) and (\ref{2.4f3})  satisfy
\begin{equation}
   -\frac{M}{t}< O\( {1 \over t} \)<\frac{M}{t}.
   \end{equation}
   This is easy to see since the  $O\( {1 /t} \)$  terms, in addition to depending on $t$, depend on $a$, $b$, $p/j$   and $p/j\le p/k\le 1$.

Using  (\ref{2.4f2}), (\ref{2.4f3})  we have  
   \begin{eqnarray}
       A_{j,k,p}(t)  \label{2.4f5}& =&{\ga \over d^{2}}\(-(d+2)   +{b(dt-b)p\over k}+{a(dt-a)p\over j} -{p^{2}\over jk}\) \label{5.15} \\
   &&\qquad +{1 \over d} \( a+b-{  p(dt-a)\over j}-{ p(dt-b)\over k}\) +O\( {1 \over t} \) \nn\\
     & =&{-(d+2)\ga+  d(a+b) \over d^{2}}\nn\\
   &&\qquad+{ (\ga a-d) p(dt-a)\over jd^{2}}+{ (\ga b-d) p(dt-b)\over kd^{2}} -{\ga p^{2}\over jkd^{2}} +O\( {1 \over t} \)  \nn\\
   & =&{-2\ze \over d^{2}}
  +{ (a-1)^{2} p(dt-a)\over jd^{2}}+{ (b-1)^{2} p(dt-b)\over kd^{2}}  -{\ga p^{2}\over jkd^{2}} +O\( {1 \over t} \).  \nn
   \end{eqnarray}
where, for the final equality we use  (\ref{2.28}).

   Note that
  \begin{eqnarray}
  B_{j,k,p}(t)&:= &{ (a-1)^{2} p(dt-a)\over j}+{ (b-1)^{2} p(dt-b)\over k}\label{2.4f6}\\
  & \geq &2p |(a-1)(b-1)|{\sqrt{(dt-a)(dt-b)} \over \sqrt{jk} }\nn\\
   & = &2p |d+2-(a+b)|{\sqrt{(dt-a)(dt-b)} \over \sqrt{jk} }\nn.
   \end{eqnarray}
   (For the inequality use $\al^{2}+\bb^{2}\ge 2\al\bb$.)
Therefore   since $ \ze =a+b-(d+2)>0$,   
 \bea
    A_{j,k,p}(t)&\ge&\frac{2}{d^{2}}\(p{\sqrt{(dt-a)(dt-b)} \over \sqrt{jk} }-1\)\ze-{\ga p^{2}\over d^{2}jk} +O\(\frac{1}{t}\)\nn\\
 &=& \frac{2}{d^{2}}\({p(dt+O(1)) \over \sqrt{jk} }-1\)\ze-{\ga p^{2}\over d^{2}jk} +O\(\frac{1}{t}\)  \label{5.18}\\
 &=&\frac{2}{d^{2}}\({pdt \over \sqrt{jk} }\(1+O\(\frac{1}{t}\)-\frac{\ga p}{ \ze\sqrt{jk}d t}\)-1\)\ze  +O\(\frac{1}{t}\)\nn.
   \eea

 Thus we see that there exists a function  $\ep_{t}$, depending only on $a$ and $b$ such that 
 \begin{equation}
      A_{j,k,p}(t) \ge\frac{2}{d^{2}}\({pdt \over \sqrt{jk} }(1-\ep_{t})-1\)\ze +O\(\frac{1}{t}\)  ,\label{5.19}
   \end{equation}
 where
 \begin{equation}
  \lim_{t\to \ff} \ep_{t}=0,\label{5.20}
   \end{equation}  
   and, as we point out above the $O\( 1/t\)  $ is independent of $p$, $j$ and $k$.
  
\begin{remark}\label{rem-5.1} {\rm We interrupt this proof to make some comments which  may be helpful in understanding what is going on. Note that 
if  
\be
{\sqrt{jk}\over dt}\le 1-\wt\ep\label{2.39}\qquad\mbox{for some $\wt\ep>0$}\label{5.21}
\ee
then  
\be
R_{j,k}(t,c)\ge R_{j,k,0}(t,c) \ge (1-\de){1-2c^2(\ze/d)\over d^{2}t^2} \al^{j-1}\bb^{k-1}\qquad\mbox{as $t\to\ff$}\label{5.22}
\ee
for all $\de>0$. This follows from (\ref{2.4fx}) and (\ref{5.19}) since when (\ref{5.21}) holds
\be 
      A_{j,k,p}(t)\ge \frac{2}{d^{2}}\( {1-\ep_{t}\over(1-\wt\ep)}-1\)\ze +O\(\frac{1}{t}\) >0 ,\label{5.19aa}
   \ee 
for all $p\ge 1$, for all $t$ is sufficiently large. Consequently when (\ref{5.21}) holds   $R_{j,k}(t,c)>0$ for all $t$ is sufficiently large when
\begin{equation} 
   c^{2}<\frac{d}{2\ze}.
   \end{equation}
(Here we also use (\ref{5.20}).)

(When $\ze\le 0$,   (\ref{5.22}) shows that $R_{j,k}(t,c)>0$ for all $c\in R^{1}$. This is   what we expect. (See the paragraph containing (\ref{4.5}).)

 }\end{remark}

 We use the next two lemmas to complete the proof of  Theorem \ref{theo-cp}, in case 1.

 \begin{lemma} \label{lem-5.2}For any  $N_{0}\in R^{+}$, we can find $c_{0}>0$ and $t_{c_{0}}<\ff$ such that for all $t\geq t_{c_{0}}$  
 \begin{equation}
 R_{j,k,p}(t,c)>0     \label{5.25}
 \ee
  for all $|c|\leq c_{0}$ and all $p$, $j$ and $k$ for which   $ jk/ t\le N_{0}$.
 \end{lemma}

\Proof   This follows from (\ref{2.4fx}) when $p=0$. Therefore, we can take $p\ge 1$. 

We first show that for any   $N\in R^{+}$,  $R_{j,k,p}(t)>0$ when $ \sqrt{jk}=Nd t$,    for all $t$ sufficiently large.
 By Remark \ref{rem-5.1} we can assume 
    that    $ N\ge 1-\wt\ep $.      
    It follows from (\ref{5.19}) that  

 \begin{equation}
      A_{j,k,p}(t) \ge \frac{2}{d^{2}}\({p \over N }(1-\ep_{t})-1\)\ze +O\(\frac{1}{t}\),\label{5.19aaxx}
   \end{equation}
   where $\ep_{t}$ satisfies (\ref{5.20}). 
 Therefore when $p\geq \La N$ for any $\La>1$,   $A_{j,k,p}(t)>0$, and hence $R_{j,k,p}(t,c)>0$,  for all $t$ sufficiently large.  
 
Now suppose that  
 \begin{equation}
 p<\La N.\label{97}
 \end{equation}
 Since $ \sqrt{jk}=N dt$ we see that  
 \begin{equation}
  {  \ga p^{2}\over d^{2} jk} \leq {\ga \La^{2} \over d^{4}t^{2} }=O(1/t^{2}),\label{5.27}
   \end{equation}
  where the $O(1/t)$ term is independent of $p$, $j$ and $k$. 
   
 Note that by   (\ref{lem-2.2})  and (\ref{5.11}) 
 \begin{equation}
 {1 \over \wt t^2} F_{j,k,p}(t)  =   {(j-p)(k-p)\over d^{2}t^{2}( p+1)jk}+   O\(\frac{1}{t^{3}}\)\label{5.1m1}
   \end{equation} 
  Therefore, 
if in addition  to (\ref{97}) we also have $\La N \le j/2 $, so that $p<j/2$, we see by (\ref{5.10}),     (\ref{5.15}) and (\ref{5.27}) that  
 \bea
 &&{1 \over \wt t^2}F_{j,k,p}(t)+{c^2\over \bar t^2} A_{j,k,p}(t)\label{3.49j}\\
&&\qquad\ge \frac{1}{d^{2}t^{2}}\( {(j-p)(k-p)\over (p+1)jk}- {c^{2}2\ze\over d}\)+O\(\frac{1}{t^{3}}\)\nn\\   
&&\qquad
   \ge\frac{1}{d^{2}t^{2}}\(  {1\over 4(p+1)}- {c^{2}2\ze\over d }\)+
   O\(\frac{1}{t^{3}}\)\nn \\
&&\qquad
   \ge\frac{1}{d^{2}t^{2}}\(  \frac{1}{8\La N }- {c^{2}2\ze\over d } \) +O\(\frac{1}{t^{3}}\).
\nn
 \eea 
Therefore we can obtain $R_{j,k,p}(t)\ge 0$ by taking 
   \begin{equation}
   c^{2}\le\frac{d }{16\La' N \ze}\label{3.51a}
   \end{equation}
 for some $\La' >\La$.

   Now suppose that  $\La N> j/2 $. In this case we use (\ref{5.15})     to see that 
 \begin{equation}
   A_{j,k,p}(t)\ge -{2 \ze\over d^{2}}+{(a-1)^2p \,t\over 2d j}+O(1/t).\label{5.31}
   \end{equation}  
   It is easy to see that the right-hand side of (\ref{5.31}) is greater than zero for all $t$ sufficiently large since 
   \begin{equation}
  {(a-1)^2p\,t\over 2 dj}\ge \frac{1}{4dN\La}(a-1)^2 t.
   \end{equation}
Thus we see that for any fixed $N$,   $R_ {j,k,p}(t) > 0$ for all $ t$  sufficiently large.  

Since the $O(1/t)$ terms are independent of $p$, $j$ and $k$ this analysis works for all $j$ and $k$ satisfying (\ref{5.25}), and all $1\le p\le j$  as long as    
(\ref{3.51a})      holds with $N$     replaced by $N_{0}$.      \qed   

\begin{lemma} \label{lem-5.3} For   all $N_{0}$ and $B\in R^{+} $      we can find   a    $c'_{0}>0$ and $t_{c'_{0}}<\ff$ such that for all $t\geq t_{c'_{0}}$  
 \be 
 R_{j,k,p}(t,c)>0 
 \ee
    for all $|c|\leq c'_{0}$ and all   $0\le p\le j\le k$ for which   
\be 
N_{0}t\leq\sqrt{jk}\leq Bt\log t.\label{5.25a} 
\ee
\el

 (The value of $N_{0}$ in Lemmas \ref{lem-5.2} and \ref{lem-5.3} can be taken as we wish. It will be assigned in the proof of this lemma.)
 
  \medskip	
\Proof    
  By adjusting  $N_{0}$ and $B$ we can replace (\ref{5.25a}) by the condition 
\begin{equation}
N_{0}\wt t\leq\sqrt{jk}\leq B\wt t\log \wt t.\label{5.25ax}
\end{equation}
   Using (\ref{2.4f5}), we see that if $  j\leq \rho \wt t$  
    \begin{eqnarray}
    A_{j,k,p}(t)&=&{-2\ze\over d^{2}}+{ (a-1)^{2} p(dt-a)\over jd^{2}}+{ (b-1)^{2} p(dt-b)\over kd^{2}} 
   \nn\\
   && \hspace{2  in} -{\ga p^{2}\over jkd^{2}} +O\( {1 \over t} \)  \label{2.4f5a}\\
   & \geq &{-2\ze\over d^{2}}+{ (a-1)^{2} p(dt-a)\over jd^{2}}- {\ga\over d^{2}}  +O\( {1 \over t} \)\nn\\
     & \geq &{-2\ze\over d^{2}}+{ (a-1)^{2} \over2 \rho d }- {\ga\over d^{2}}  +O\( {1 \over t} \)\nn.
      \end{eqnarray}    
   Clearly, there exists a $\rho>0$, independent of $j$ and $k$ such that this term  is positive.  Thus we can assume that
    \begin{equation}
   j\geq \rho \wt t.\label{2.40w}
   \end{equation}

   Furthermore, when  $  \sqrt{jk}/\wt t= N$,  it follows from  (\ref{2.51a1}) that we can write (\ref{5.18}) as  
  \be  A_{j,k,p}(t)\ge \frac{2}{d^{2}}\({p  \over N}\(1 +O\(\frac{1}{t}\)-\frac{\ga  }{   \ze (dt)^{2}}\)-1\)\ze +O\(\frac{1}{t}\).\label{5.37}
 \ee  
 Let $\de_{N}=(10\log N/N)^{1/2}$. Clearly, if $p>(1+\de_{N})N$, the right-hand side of (\ref{5.37}) is positive for all $t$ sufficiently large.
   Therefore, when   $  \sqrt{jk}/\wt t=N$, we may assume that   
          \begin{equation}
p\leq (1+\de_{N})N\label{2.40q}.
   \end{equation}
 (The value  chosen for $\de_{N}$  simplifies calculations made later in this proof.)
   
   In addition  we can also assume that 
   \be
    p\geq p_{0}
    \ee
 for any finite $ p_{0}$, since if $p< p_{0}$
\bea 
  F_{j,k,p}(t)\ge  F_{j,k,p_{0}}(t)\ge c ^{2}A_{j,k,p}(t).\label{5.39}
 \eea
 for all $c>0$ sufficiently small.  
 
 \medskip	
 We use the next lemma in the proof of Lemma  \ref{lem-5.3}.

   \bl\label{lem-3.9}   For $j\le k$, with  $p$ and $j$  large and $p/j$ small  
 \bea
 &&
 {j\choose p}{k\choose p}= { 1\over  2\pi p}{\({e^2jk\over  p^2}\)^{p}}\label{2.40}\\&&\qquad\exp\(-{p\over 2j}(p-1)-{p\over 2k}(p-1)+O(p^3/j^2)\)\(1+O(p^{-1})\) .\nn
 \eea 
 When  $\wt t\in R^{+}$ is large and   $\sqrt{jk}/\wt t=N$, under  assumptions (\ref{2.40w}) and (\ref{2.40q})
 \bea
 &&
 {1\over \wt t^{2p}}{j\choose p}{k\choose p}= {1 \over 2\pi p}\({eN\over p}\)^{2p}\(1+O(p^{-1})\).\label{2.51}
 \eea
  
 \el

 \Proof By Stirlings's formula  for integers   $q$,
 \be
 q!=\sqrt{2\pi}q^{q+1/2}e^{-q}\(1+O(q^{-1})\).\label{5.47}
 \ee
 Therefore, since $j$ is large and $p/j$ is small, terms of the form 
 \begin{equation}
{\(1+O(j^{-1})\) \over  \(1+O(p^{-1})\)\(1+O((j-p)^{-1})\)}=\(1+O(p^{-1})\)\label{2.40n}.
 \end{equation}
Using this we see that  
  \bea
  {j\choose p}&&={j! \over (j-p)!p!}\\
  & &={1 \over \sqrt{2\pi}}{j^{j+1/2} \over (j-p)^{(j-p+1/2)}\,p^{p+1/2}}\(1+O(p^{-1})\) \nn\\
    & &={1 \over \sqrt{2\pi\,p}}\(\frac{ j}p\)^p{j^{j-p+1/2} \over (j-p)^{(j-p+1/2)}}\(1+O(p^{-1})\) \nn\\
        & &={1 \over \sqrt{2\pi\,p}}\(\frac{ j}p\)^p{1\over \(1-{p \over j}\)^{(j-p+1/2)}}\(1+O(p^{-1})\) \nn\\
& &={1 \over \sqrt{2\pi\,p}}\(\frac jp\)^p e^{-(j-p+1/2)\log(1-p/j)}\(1+O(p^{-1})\) \nn\\
 & &={1 \over \sqrt{2\pi\,p}}\(\frac jp\)^p e^{(j-p+1/2)(p/j+p^2/(2j^2) +O(p^3/j^3))}\(1+O(p^{-1})\) \nn\\
 & &={1 \over \sqrt{2\pi\,p}}\(\frac{e j}p\)^p e^{( -p^2/(2j )+p/(2j)
 +O(p^3/j^2))}\(1+O(p^{-1})\) \nn\\
 & &={1 \over \sqrt{2\pi\,p}}\(\frac{e j}p\)^p e^{( -p(p-1)/(2j )
 +O(p^3/j^2))}\(1+O(p^{-1})\) \nn.
 \eea
Since this also holds with $j$ replaced by $k$ we get (\ref{2.40}). 

To get (\ref{2.51}) we multiply each side of   (\ref{2.40}) by $\wt t^{-2p}$ and  substitute for $\sqrt{jk}/\wt t=d N$    and use  the fact that under the assumptions (\ref{2.40w}) and (\ref{2.40q}), 
 \begin{equation} 
 {p^{3} \over j^{2}}\leq  { p^{2} \over j }\leq   \frac{(1+\de_{N})^{2}N^{2}}{\rho t} \label{6.127}.
 \end{equation} 
Consequently, for all $t$ sufficiently large
 \begin{equation}
 \exp \(-{p\over 2j}(p-1)-{p\over 2k}(p-1)+O(p^3/j^2) \)=1+ O\(N^{2}/t\).\label{6.128}
 \end{equation}
 \qed

 \medskip	
\noindent 
{\bf Proof of Lemma \ref{lem-5.3}   continued} We show that under the assumptions  (\ref{2.40w}) and (\ref{2.40q}),  when $\sqrt{jk}/  \wt t=N$,    for   
 $N_{0}\leq N\leq B\log   \wt t$,  for any $0<B<\ff$, and $t$ is sufficiently large,
 \begin{equation}
 { \wt  t^{2}\over   \al^{j-1}\bb^{k-1}}\sum_{p=p_{0}}^{(1+\de_{N})N}P_{j,k,p}(t)\geq Ce^{2N}  {1\over N^{3/2} } \label{6.11}
 \end{equation}
 for some $C>0$, independent of $N$,
 and 
  \begin{equation}
 {\bar t^2\over    \al^{j-1}\bb^{k-1}}\sum_{p=p_{0}}^{(1+\de_{N})N}Q_{j,k,p}(t)\geq -De^{2N}  {1\over N^{3/2} } \label{6.12}
 \end{equation}
 for some $D<\ff$, independent of $N$. If  (\ref{6.11}) and (\ref{6.12}) hold, we can find a $c_{0}>0$ such that for all  $c^{2}\le c_{0}^{2}$ 
   \begin{equation}
 \sum_{p=p_{0}}^{(1+\de_{N})N}R_{j,k,p}(t,c)\geq 0\label{6.13}
 \end{equation}
 for all  $t$ sufficiently large.  Since we have already established that $A_{j,k,p}(t)>0$,  when $p<p_{0}$ and $p>(1+\de_{N})N$,   this  completes the proof of  Lemma \ref{lem-5.3}. 
 
  Thus it only remains to prove (\ref{6.11}) and (\ref{6.12}). We do (\ref{6.11}) first. It is considerably easier than  (\ref{6.12}).  By  Lemma \ref{lem-P} and  (\ref{2.51})
  \bea
 &&{\wt t^{2}\over   \al^{j-1}\bb^{k-1}}\sum_{p=p_{0}}^{(1+\de_{N})N}P_{j,k,p}(t)\label{6.14}\\
 &&\qquad=\sum_{p=p_{0}}^{(1+\de_{N})N}
 {1\over \wt t^{2p}}{j\choose p}{k\choose p} {(j-p)(k-p)\over (p+1)jk}\nn\\
 &&\qquad\ge C \sum_{p=p_{0}}^{(1+\de_{N})N}
 {1\over   \wt t^{2p}}{j\choose p}{k\choose p} {1\over p} \nn\\
 &&\qquad\ge C \sum_{p=p_{0}}^{(1+\de_{N})N}
 {1\over p^{2}}\({eN\over  p  }\)^{2p  } \nn.
 \eea
 
In order to calculate this last sum we   consider the function  

\begin{equation}
f_{m}(y)={1 \over y^{m}}\({eN\over  y  }\)^{2y  }={1 \over y^{m}}e^{2y(1+\log N-\log y)}\label{jr.10}
\end{equation}
for $m\geq 0$ and $y\geq 2.$
We have    
\bea
f_{m}'(y)&=&\(  {-m \over y}   +2(1+\log N-\log y)-2   \)f(y) \label{jr.11}\\
&=&\(  {-m \over y}   +2(\log N-\log y)  \)f(y).\nn
\eea
This has a unique root $y_{m}$ where
\begin{equation}
\log y_{m}+{m \over 2y_{m}}=\log N.\label{jr.12}
\end{equation}
(Clearly, $y_{0}=N$).  
Let $y_{m}=N(1+\ep_{m})$.  Then\begin{equation}
\log (1+\ep_{m})+{m \over 2N(1+\ep_{m})}=0.\label{jr.14}
\end{equation}
Consequently
\begin{equation}
\ep_{m}=-{m \over 2N}+O(N^{-2}),\label{jr.15}
\end{equation}
which implies that   
\begin{equation}
y_{m}=N-{m \over 2}+O(1/N).   
   \end{equation}
Making use of the fact that $\(  {-m /y_{m}}   +2(1+\log N-\log y_{m})-2   \)=0$, we see that  
\bea
f_{m}''(y_{m})&=&\(  {m \over y_{m}^{2}}   -\frac{2}{y_{m}}  \)f(y_{m}) <0.\label{jr.11x}
\eea
Therefore   
\begin{equation}
\sup_{y\geq  2}f_{m}(y ) =f_{m}(y_{m})\leq  {1\over (N-m)^{m}}e^{2N}.\label{jr.16}
\end{equation}
We also note that since $ {-m / y} -2\log y$ is increasing for $y>m/2$,   $f_{m}'(y)$ is positive  for $m/2<y<y_{m}$ and negative for $y>y_{m}$. Consequently, $f_{m}(y) $ is unimodal.

\medskip	
Now consider the last line of (\ref{6.14}). The function being summed is $f_{2}(p)$. The above discussion shows that this function is unimodal, with a maximum at, at most,  two points at which it is less than $ 2 e^{2N}/ N^{2}$.  Consequently, to obtain (\ref{6.11}) we can replace the sum in the last line of (\ref{6.14}) by an integral and show that
 \begin{equation}
 I_{1}:=\int_{ p_{0}}^{(1+\de_{N})N}
 {1\over r^{2}}\({eN\over  r  }\)^{2r  }\,dr\geq Ce^{2N}  {1\over N^{3/2} } \label{6.11a}.
 \end{equation}
Making the change of variables $r=xN$  we have
  \bea
 I_{1} ={  1\over N} \int_{ p_{0}/N}^{1+\de_{N}}
 {1\over x^{2}}\({e\over x  }\)^{2xN  }\,dx.   \nn
 \eea
 Recall  that $N_{0}\leq N\leq 2\log \wt t$, and  that we can take $N_{0}$ as large as we want, (but fixed and independent of $t$), and that $\de_{N}=(10\log N/N)^{1/2}$. Therefore
 \bea
 I_{1}&\geq  &
 {1\over N}\int_{1-(10\log N/N)^{1/2}}^{1+(10\log N/N)^{1/2}}  { 1\over  x^{2}}{\({e \over x}\)^{2xN}}\,dx\label{2.51h} \\
 &\geq  &
 {1 \over 2 N}\int_{1-(10\log N/N)^{1/2}}^{1+(10\log N/N)^{1/2}} {\({e \over x}\)^{2xN}}\,dx. \nn
 \eea

 We write 
  \be
  \({e\over x}\)^{2xN}=\exp\(2xN(1-\log x)\)\label{j.58}.
  \ee
  Set $x=1+y$ and
     note that for $|y|$ small
      \begin{eqnarray} 
     x(1-\log x) \label{j.58j} 
  &=& (1+y)(1-\log (1+y) )\\
  & = &(1+y)\(1-\sum_{n=1}^{\ff}(-1)^{n-1}{y^{n} \over n}\)  \nonumber\\
  & = &1+y -\sum_{n=1}^{\ff}(-1)^{n-1}{(1+y)y^{n} \over n} \nonumber\\
  & = &1 -\sum_{n=2}^{\ff}(-1)^{n-1}y^{n}\({1 \over n}-{1 \over n-1}\) \nonumber\\
  & = &1 -\sum_{n=2}^{\ff}(-1)^{n} {y^{n} \over n(n-1)}.  \nonumber
  \end{eqnarray}

  When $|y|\leq (10\log N/N)^{1/2}$,  so that $|y|^{3}N <<1$, this shows that
    \bea
  \({e\over x}\)^{2x N} 
  &=&e^{2 N}  e^{ -y^2 N +O(y^3N)}\label{2.59}\\
  & = &   e^{2 N}  e^{ -y^2 N }\(1+O(|y|^{3}N )\).\nn
 \eea 
It follows from this that when we make the change of variables $x=1+y$ in (\ref{2.51h}) we get
  \bea
 && I_{1} \geq 
 {  e^{2 N}\over 2 N}\int_{-(10\log N/N)^{1/2}}^{(10\log N/N)^{1/2}} e^{ -y^2 N }\,dy  \\
&&\qquad \geq 
 { e^{2 N}\over 2\sqrt{2} N^{3/2}}\int_{-(20\log N)^{1/2}}^{(20\log N)^{1/2}} e^{ -y^2/2  }\,dy.\nn
 \eea
Since 
\begin{equation}
  \int_{(20\log N )^{1/2}}^{\ff} e^{-u^{2}/2}\,du \le N^{-10}\label{2.51j},
\end{equation} 
 we see that (\ref{6.11a}) follows.   Thus we have established (\ref{6.11}).

\medskip	
Before proceeding to the proof of (\ref{6.12}) we note that 
\begin{equation}
  \sum_{p=p_{0}}^{(1+\de_{N})N}\label{5.72}
 {1\over \wt t^{2p}}{j\choose p}{k\choose p}\le   { e^{2 N} \over 2 N^{1/2}}.
   \end{equation}
 To prove this we use  (\ref{2.51}) and the same argument that enables us to move from a sum to an integral that is given in (\ref{6.14})--(\ref{6.11a}), except that we use (\ref{jr.16}) with $m=1$. We continue  and then  use (\ref{2.59}) to get
\begin{eqnarray}
&&\sum_{p=p_{0}}^{(1+\de_{N})N}
 {1\over \wt t^{2p}}{j\choose p}{k\choose p}
\label{jr.7}\\
 &&\qquad\leq  \int_{p_{0}}^{(1+\de_{N})N}{1 \over u} \({eN\over  u }\)^{2u}\,du+O\({e^{2 N} \over N}\)\nn\\
 &&\qquad\leq  e^{2 N}  \int_{0}^{\de_{N}}   e^{ -y^2 N }\,dx+O\({e^{2 N} \over N}\)\leq   { e^{2 N} \over 2N^{1/2}} . \nn
\end{eqnarray}

 \medskip	
 We now obtain (\ref{6.12}).  When $\sqrt{jk}=\wt t N$, by (\ref{5.10}) and   (\ref{5.37}),
   \bea
 &&{\bar t^2\over    \al^{j-1}\bb^{k-1}}Q_{j,k,p}(t)\\
 &&\qquad={1\over \wt t^{2p}}{j\choose p}{k\choose p}A_{j,k,p}(t)\label{6.120w}\nn\\
  & &\qquad\ge{2\over  d^{2}\wt t^{2p}}{j\choose p}{k\choose p}\lc 2\({p  \over N }-1\)\ze +O\(\frac{1}{t}\)\rc \nn.
 \eea

By  (\ref{jr.7}) we see that
\begin{equation}
   \sum_{p=p_{0}}^{(1+\de_{N})N}{1\over \wt t^{2p}}{j\choose p}{k\choose p}O\(\frac{1}{t}  \)=O\({e^{2 N}\over t}\)
   \end{equation}
Therefore,  to obtain (\ref{6.12}), it suffices to show that for some $D<\ff$
\begin{equation}
 \sum_{p=p_{0}}^{(1+\de_{N})N}{1\over \wt t^{2p}}{j\choose p}{k\choose p}  \({p  \over N }-1\) \geq -De^{2N}  {1\over N^{3/2} }\label{6.12a}.
\end{equation}
  Here we use the fact that $N\le B\log \wt t$ for some $0<B<\ff$.

\br {\rm
 Since the proof of  (\ref{6.12a}) is rather delicate we make some  heuristic comments to explain how we proceed. When   $\sqrt{jk}=\wt t N$ the term $\wt t^{-2p}{j\choose p}{k\choose p}$, as a function of $p$, is $\exp(2N)$ times values that are sort of normally distributed with mean $p=N$, and, roughly speaking,  
 \begin{equation}
   \sum_{p=p_{0}}^{(1+\de_{N})N}{1\over \wt t^{2p}}{j\choose p}{k\choose p}\sim Ce^{2N}\frac{1}{N^{1/2}},\label{jr.22}
   \end{equation}
 for all $t$ sufFiciently large. (In fact the upper bound is given in (\ref{jr.7}).)
 This is too large to enable us to get (\ref{6.12a}) so we must make use of the factors $ \({p  \over N }-1\)$, which   is an odd function with respect to  $p=N$, to get the cancellations that allow us to obtain  (\ref{6.12a}).  However, because we are canceling terms, we must take account of the error in Stirling's approximation; (see (\ref{5.47})). To do this we need to show that the estimate in (\ref{6.12a}) remains the same even when we eliminate the terms in the summand that are not close to $N$.}
 \er
 
 \medskip	
\noindent
{\bf Proof of Lemma \ref{lem-5.3}  continued }
Note that by  (\ref{2.51})
\be
\sum_{p=p_{0}}^{N(1-N^{-1/4})}{1\over \wt t^{2p}}{j\choose p}{k\choose p}   
\label{jr.21} \leq C  \sum_{p=p_{0}}^{N(1-N^{-1/4})}f_{1}(p).  \label{5.79}
\ee
The fact that $f_{m}(y)$ is unimodal on $y>m/2$ implies that $f_{1}(p)$ is increasing on the interval $[p_{0}, N(1-N^{-1/4})]$. Therefore   
\begin{eqnarray}
\sum_{p=p_{0}}^{N(1-N^{-1/4})}{1 \over p} \({eN\over  p }\)^{2p}   
&\leq &  C  N f_{1}(N(1-N^{-1/4}))\label{jr.21xx}\\
&\leq & C \({e\over  1-N^{-1/4} }\)^{2N(1-N^{-1/4})}. \nonumber\\
&= & C e^{2N(1-N^{-1/4})(1-\log (1-N^{-1/4} ))}. \nonumber\\
&\leq  & C e^{2N(1-N^{-1/4})(1+N^{-1/4})}=Ce^{2N -2N^{1/2} }. \nonumber
\end{eqnarray}

Let $\de'_{N}=N^{-1/4}$. The argument immediately above shows that to prove (\ref{6.12a}), it suffices to show that
\begin{equation}
J_{1}:=\sum_{p=(1-\de'_{N})N}^{(1+\de_{N})N}{1\over \wt t^{2p}}{j\choose p}{k\choose p}  \({p  \over N }-1\) \geq -De^{2N}  {1\over N^{3/2} }.\label{6.12b}
\end{equation}

\medskip	 By (\ref{2.51})
\begin{equation}
J_{1}=\frac{1}{2 \pi }\sum_{p=(1-\de'_{N})N}^{(1+\de_{N})N} \({eN\over   p }\)^{2p}  \({1 \over N }-\frac{1}{p}\)\(1+O\(\frac{1}{p}\)\).\label{6.12c}
\end{equation}
Using  (\ref{jr.7}) together with the fact that since $p\geq (1-\de'_{N})N$,  
$1/p \leq  1/N $, we see that
\bea
 \sum_{p=(1-\de'_{N})N}^{(1+\de_{N})N} \({eN\over   p }\)^{2p}  \({1 \over N }-\frac{1}{p}\) O\(\frac{1}{p}\)&\ge&- \sum_{p=(1-\de'_{N})N}^{(1+\de_{N})N} \({eN\over   p }\)^{2p}  \frac{1}{p} \, O\(\frac{1}{p}\)\nn\\
 &\ge&-C   e^{2N}{1  \over N^{3/2}}  .\label{6.12df}\nn
\eea
Therefore, to obtain (\ref{6.12b})
 that it suffices to show that
\begin{equation}
 \sum_{p=(1-\de'_{N})N}^{(1+\de_{N})N} \({eN\over   p }\)^{2p}  \({1 \over N }-\frac{1}{p}\) \geq -De^{2N}  {1\over N^{3/2} }.\label{6.12d}
\end{equation}

\medskip

In a minor modification of the analysis of $f_{m}(y)$, we write 
\begin{equation}
 h(y):=   \({eN\over   y}\)^{2y}  \({1 \over N }-\frac{1}{y}\) =\exp\(2y\(1+\log N- \log y\)\)\({1 \over N }-\frac{1}{y}\).
   \end{equation}
 Therefore
 \begin{equation}
   h'(y)=\(\(2 \(1+\log N-\log y\) -2 \)\({1 \over N }-\frac{1}{y}\)+\frac{1}{y^{2}} \)\({eN\over   y }\)^{2y}. 
   \end{equation}
Let $y=(1+\om)N$. Then $ h'(y)=0$ when
\begin{equation}
  - 2\om  \log(1+\om)   +\frac{1}{N (1+\om) }=0.
   \end{equation}
 This equation is satisfied when 
  \begin{equation}
  \om=\pm \frac{1}{\sqrt{2N} }+O\(\frac{1}{N}\).
   \end{equation}
 Note that when $y=(1+\om)N$
 \begin{equation}
    \({eN\over   y }\)^{ y}= \({e \over  1+\om}\)^{  (1+\om)N}\le e^{N},
   \end{equation}
 because $(e/x)^{x}$ is maximized when $x=1$. Therefore
 \begin{equation}
    \({eN\over   y }\)^{2y}  \({1 \over N }-\frac{1}{y}\) \le  {e^{2N}\over N}\frac{\om}{1+\om}
   \end{equation}
 from which we get  
 \begin{equation}
 \sup_{1\leq y\leq (1+\de_{N})N}|h(y)|\leq  C\(\frac{e^{2N}}{N^{3/2}}\).\label{5.68}
 \end{equation}
 
It is easy to see that $h(y)$ is negative for $1\le y\le N$ and that it decreases to its minimum value at $N(1-\om)$ and then increases to zero at $y=N$. It then increases to its maximum value at $N(1+\om)$ and then decreases for $N(1+\om)\le y\le (1+\de_{N})N$.
Consequently the difference between
\begin{equation}
   \sum_{p=(1-\de'_{N})N}^{{(1+\de_{N})N}}h(p) \qquad\mbox{and}\qquad    \int_ {(1-\de'_{N})N} ^{{(1+\de_{N})N}}h(p)\,dp
   \end{equation}
 differs by at   most $4\max_{1\le p\le (1+\de_{N})N}|h(p)|$. Since this is $O(e^{2N}/N^{3/2})$ by (\ref{5.68}), and we are only trying to obtain (\ref{6.12d}), we can neglect this discrepancy.  Therefore to obtain (\ref{6.12d}) we need only show that
 \begin{equation}
  \int_ {(1-\de'_{N})N} ^{{(1+\de_{N})N}}\frac{1}{p} \({eN\over   p }\)^{2p}  \({p \over N }-1\) \,dp\ge -{D'}\frac{e^{2N}}{N^{3/2}}.\label{5.76}
   \end{equation}
 Under the change of variables $p=xN$ the integral in (\ref{5.76}) is equal to
 \begin{equation}
     \int_ {1-\de'_{N}} ^{{1+\de_{N}}}\frac{1}{x} \({e \over  x}\)^{2xN}  \(x-1 \) \,dx .\label{5.77}
 \end{equation}

 By (\ref{j.58j}),  in which $x=1+y$  ,
    \bea
  \({e\over x}\)^{2x N} 
  &= &e^{2 N}  e^{ -y^2 N +y^3 N/3 +O(y^4) N}\label{2.59sas}\\
  & =&  e^{2 N}  e^{ -y^2 N }\(1+ {y^3 N\over 3}+O(y^4) N\)\nn.
 \eea 
 Therefore, with the change of variables $x=1+y$ we write the integral in (\ref{5.77}) as
 \be 
      e^{2 N}\int_{-\de'_{N}}^{\de_{N}} \frac{y}{1+y  }e^{ -y^2 N } \(1+ {y^3 N\over 3}+O(y^4) N\)\label{5.82}
 \,dy. 
   \ee 
   We use $(1+y)^{-1 }=(1-y +y^{2} -y^{3} +O(y^{4}) )$  to write
  \bea
    &&\frac{y}{1+y }  \(1+ {y^3 N\over 3}+O(y^4) N\)\\
    &&\qquad\quad=y- y^{2} +  y^{3} +\frac{y^{4 }N}{3}- y^{4 }+O(y^{5})N.\nn
   \eea
   Using this we see that   (\ref{5.82})
   \begin{equation}
  =    e^{2 N}\int_{-\de'_{N}}^{\de_{N}} e^{ -y^2 N } \(y- y^{2} +  y^{3} +\frac{y^{4 }N}{3}- y^{4 }+O(y^{5})N\)\label{5.82q}
 \,dy. 
   \end{equation}
   Recall that   $\de_{N}=(10\log N/N)^{1/2}$ and $\de'_{N}=N^{-1/4}$ . Since
    \begin{equation}
  e^{2 N} \int_{(10\log N/N)^{1/2}}^{\ff}e^{-y^{2}N}\,dy  \le {e^{2 N}\over N^{10}}\label{5.98}
   \end{equation}
   and
       \begin{equation}
  e^{2 N} \int_{-\ff}^{-N^{-1/4}}e^{-y^{2}N}\,dy  \le e^{ 2N-N^{1/2}} ,
   \end{equation}
  errors we can ignore in obtaining (\ref{6.12}), we can simplify matters by replacing the integral in (\ref{5.82q}) by 
      \bea
&&    e^{2 N}\int_{-\ff}^{\ff} e^{ -y^2 N } \(y- y^{2} +  y^{3} +\frac{y^{4 }N}{3}- y^{4 } +O(y^{5})N\)\label{5.82qww}
 \,dy\\
 &&\qquad   =- e^{2 N}\int_{-\ff}^{\ff} e^{ -y^2 N } \(  y^{2}-\frac{y^{4 }N}{3}+y^{4 } +O(y^{5})N\) 
 \,dy\nn\\
 &&\qquad   =- e^{2 N}\int_{-\ff}^{\ff} e^{ -y^2  } \( { y^{2} \over N^{3/2}}-\frac{y^{4 }}{3N^{3/2}}+{y^{4 } \over N^{5/2}} +O(y^{5})N^{-2}\) 
 \,dy\nn\\
 &&\qquad   =- {e^{2 N} \over N^{3/2}}\int_{-\ff}^{\ff} e^{ -y^2  } \( y^{2}-\frac{y^{4 }}{3} \) 
 \,dy+O\({e^{2 N} \over N^{2}}\)\nn.
   \eea

Since 
\bea
&&
{1 \over \sqrt{\pi}}\int_{-\ff}^{\ff} e^{ -y^2  } \( y^{2}-\frac{y^{4 }}{3} \) 
 \,dy\label{5.76f}\\
 &&\qquad={1 \over \sqrt{2\pi}}\int_{-\ff}^{\ff} e^{ -y^2/2  } \( \frac{y^{2 }}{2}-\frac{y^{4 }}{12} \) 
 \,dy\nn\\
 &&\qquad={1 \over 2}-  \frac{1}{4}={1 \over 4}\nn,
\eea
 we obtain (\ref{5.76}).\qed
 
 \medskip	\noindent\label{sec5page}
{\bf Proof of Theorem \ref{theo-cp}  when {$\bf c_{1}=c_{2}=c$} and {$\bf  EG_{1} G_{2}>0 $} concluded }  \label{41}  
Consider the Gaussian random variable $ (G_{1}/\ga,G_{2}/\ga)$ where $\ga=(E G_{1}G_{2}) ^{1/2}$. This random variable has covariance $\Ga$ in (\ref{cov}). By Lemma \ref{lem-5.3} there exists a $c'_{0}>0$ such that $ (G_{1}/\ga+c,G_{2}/\ga+c)$ has infinitely divisible squares for all $|c|\le c'_{0}$.  Let $\wt c$ be the supremum of the $c_{0}' $ for which this holds. Since, by hypothesis, (\ref{1.12q}) does not hold, $\wt c$ is finite.  
Therefore, $ (G_{1}/\ga+c,G_{2}/\ga+c)$ has infinitely divisible squares for all $|c|<\wt c$ and not for any $c$ or which $|c|>\wt c$. Translating this into the notation used in Theorem \ref{theo-cp} we have $ (G_{1}/\ga+c\al,G_{2}/\ga+c\al)$ has infinitely divisible squares for all $|\al|<\wt c/c$ and not for any $|\al|$ for which $|\al|>\wt c/c$.

Therefore, to complete the proof of Theorem \ref{theo-cp}  when {$  c_{1}=c_{2}=c$} and {$   EG_{1} G_{2}>0 $ we need only show that $ (G_{1}/\ga+c,G_{2}/\ga+c)$ has infinitely divisible squares for $|c|=\wt c$.  Consider the Laplace transform of  $ (G_{1}/\ga+c,G_{2}/\ga+c)$ in (\ref{idpr.2}). Since it only depends on $c^{2}$ we can simplify the notation by taking   $c>0$. 
Let $c_{m}\uparrow \wt c$.
 Abbreviate the third line of (\ref{idpr.2}) by 
$  \exp\(P+c^{2}Q\)$.
Thus $  \exp\(P+c_{m}^{2}Q\)$ is the Laplace transform of an infinitely divisible random variable. Therefore,   for each $t>0$ the power series expansion of $P+c_{m}^{2}Q$ in $s_{1}$ and $s_{2}$ has positive coefficients, except for the constant term. Thus if we write
\[P=\sum_{j,k } a_{j,k} s_{1}^{j}s^{k}_{2},\hspace{.2 in}Q=\sum_{j,k } b_{j,k} s_{1}^{j}s^{k}_{2}\]
we see that $a_{j,k}+c_{m}^{2} b_{j,k}\geq 0$ for each $(j,k)\neq (0,0)$. Letting $c_{m}\uparrow \wt c$ we therefore have $a_{j,k}+\wt c^{2} b_{j,k}\geq 0$ for each $(j,k)\neq (0,0)$. This shows that $  \exp\(\(P+\wt c ^{2}Q\) \)$  is the Laplace transform of an infinitely divisible   random variable. \qed  

\begin{remark} \label{rem-5}{\rm In the remainder of this paper we continue to prove Theorem \ref{theo-cp} for all $ c_{1},c_{2} $ and arbitrary covariance $EG_{1}G_{2}$. In each case, as immediately above, because  (\ref{1.12q}) does not hold, there exists a $  c'<\ff$ such that $ (G_{1} +c c_{1},G_{2} +cc_{2})$ does not have infinitely divisible squares for all $c$ such that $|c|>  c'$. 
Therefore, if we can show that there exists some $c\ne 0$ for which both 
\be
 (G_{1} +c c_{1},G_{2} +cc_{2}) \quad\mbox{and}\quad  (G_{1} -c c_{1},G_{2} -cc_{2})\label{5.100}
 \ee
  have infinitely divisible squares, we can use the arguments in the preceding three paragraphs to show that there exists a critical point $\wt c$ such that $ (G_{1} +c c_{1},G_{2} +cc_{2})$ has  infinitely divisible squares for all  $|c|\le \wt c$ and not for $|c|> \wt c$. Consequently,  in the remainder of this paper,  in which we consider different cases of  $ c_{1},c_{2} $ and arbitrary covariance $EG_{1}G_{2}$ we will only show that (\ref{5.100}) holds for some $c\ne 0$.
 }\end{remark}

   \section{ Proof of Theorem  \ref{theo-cp}  when {$\bf (c_{1},c_{2})=(c,\pm c)$} }\label{sec-6}
   
  We first assume that  $EG_{1}G_{2}>0$ and that   $ (c_{1},c_{2})=(c,- c)$. In this case we have  
    \bea
\lefteqn{E_G\(e^{ -( \la_1(G_{1}+c   )^2+\la_2(G_2-c   )^2)/2}\)\label{idsqaspr.2qxx}}
\\
&& =\frac{1}{(H(a,b,\la_1,\la_2))^{1/2} }   \exp\(-  c^{2}\(\frac{\rho\la_{1}\la_{2}+\la_{1}+\la_{2}}{2H(a,b,\la_1,\la_2)} \) \)\nn,
\eea
where $\rho=a+b+2$. This is exactly the same as (\ref{idpr.2}) except that $\ga$ is replaced by $\rho$. We now trace the proof in Sections \ref{sec-3}--\ref{sec-5} and see what changes. Obviously much remains the same.   In particular the power series $P$ is unchanged. The basic expression for $Q$ in (\ref{3.10})
 is essentially the same except that $\ga$ is replaced by $\rho$.  Thus Lemma \ref{lem-2.1j} is also essentially the same except that $\ga$ is replaced by $\rho$.
 
 The analysis in Section \ref{sec-4} only uses the fact that $\ga<\ff$,   and  since $\rho<\ff$, Lemma \ref{lem-1.1} also holds in this case.

In going through Section \ref{sec-5} we see the coefficients of $Q$ change, but they still  lead to essentially the same inequalities that allow us to complete the proof.
In place of (\ref{2.28}) we have
\bea
-3\rho+a+b&=&-2(3+(a+b)):=-2\wt \ze\label{6.2}\\
a\rho-1&=&(a+1)^2\nn\\
b\rho-1&=&(b+1)^2\nn.
\eea
Using this in (\ref{5.4}) and (\ref{3.23}), with $\ga$  replaced by $\rho$, we get 
 \bea
Q_{j,0}(t )  &=& \( {  (b+1)^{2} \over dt}+O(1/t^{2})\)\al^{j-1},
 \eea 
and
 \bea
Q_{0,k}(t )  &=&  \( {  (a+1)^{2} \over dt}+O(1/t^{2})\)\al^{k-1} .\label{3.23xx}
 \eea 

We also see that we get (\ref{2.4f}) with $\ga$  replaced by $\rho$ and consequently, in place of (\ref{2.4fx}), we get 
   \begin{equation}
  {R_{j,k,0}(t,c) \over  \al^{j-1}\bb^{k-1}}= {1-2c^2(\wt\ze/d)+O\( t^{-1} \)\over d^{2}t^2}\label{2.4fxq}.
   \end{equation}
 Of course the key term in the proof is the analogue of   $A_{j,k,p}(t)$. We get the third line of (\ref{5.15}) with $\ga$  replaced by $\rho$, which by (\ref{6.2}) leads to (\ref{5.18}) with $\zeta$ replaced by $\wt\zeta$ and $\ga$  replaced by $\rho$. 
Therefore,   all the subsequent lower bounds for $A_{j,k,p}(t)$ that are in Section \ref{sec-5} hold   when $\ze$ is replaced by $\wt \ze$.
In the   proof  of (\ref{6.12}) in Section \ref{sec-5} the only property of $\ze$ that is used is that is is positive. Since $\wt\zeta$ is also positive the same argument  completes the proof of Lemma \ref{lem-5.3} and   consequently, by Remark \ref{rem-5},    of Theorem  \ref{theo-cp},    when $EG_{1}G_{2}>0$ and     $ (c_{1},c_{2})=(c,- c)$. 

When $EG_{1}G_{2}<0$ and     $ (c_{1},c_{2})=(c,- c)$ we note that
\begin{equation}
   ((G_{1}+c)^{2},(G_{2}-c)^{2})\stl  ((G_{1}+c)^{2},(-G_{2}+c)^{2}).
   \end{equation}
Now $EG_{1}(-G_{2})>0$ and we are in the case proved on page \pageref{41}. Therefore, by Remark \ref{rem-5}, Theorem \ref{theo-cp} holds in this case. 

Finally when $EG_{1}G_{2}<0$ and     $ (c_{1},c_{2})=(c, c)$ we note that 
\begin{equation}
   ((G_{1}+c)^{2},(G_{2}+c)^{2})\stl  ((G_{1}+c)^{2},(-G_{2}-c)^{2}).
   \end{equation}
 Now $EG_{1}(-G_{2})>0$ and we are  in the case proved in the beginning of this section.
  \qed

   \section{ Proof of Theorem  \ref{theo-cp} when {$\bf (c_{1},c_{2})=(c,0)$}} \label{sec-7}
   
 We first assume that  $EG_{1}G_{2}>0$.  In this case we have  
    \bea
\lefteqn{ E_G\(e^{ -( \la_1(G_{1}+c  )^2+\la_2G_2^2)/2}\)\label{idsqaspr.2qyy}}
\\
&& =\frac{1}{(H(a,b,\la_1,\la_2))^{1/2} }   \exp\(-  c^{2}\(\frac{b\la_{1}\la_{2}+\la_{1} }{2H(a,b,\la_1,\la_2)} \) \)\nn.
\eea 
The term in the numerator of the exponential lacks the $\la_{2}$ that is present in  (\ref{idpr.2}) and (\ref{idsqaspr.2qxx}). Therefore, the formulas for the coefficients of the power series for the analogue of $Q$, which we denote by  $\wt Q$, are different.  It is easy to see that in place of  (\ref{as1.3})
we get \bea
\wt Q(s_1,s_2,t)
&=&-\frac{ (b t^2(1-s_1)(1-s_2)+t(1-s_1 ))}{\bar t^2}\label{as1.3xx}\\
&&\qquad
\cdot \sum_{n=0}^\ff  (\al s_1+\bb s_2-\th\al\bb s_1s_2)^n\nn\\
&=&-\frac{ (b t^2+ t-(b t^2+t)
s_1-b t^2 s_2 +b t^2s_{1} s_2 )}{\bar t^2}\nn\\ &&\qquad
\cdot\sum_{n=0}^\ff  (\al s_1+\bb s_2-\th\al\bb s_1s_2)^n\nn .
\eea
 Using this, in place of Lemma \ref{lem-2.1j}, we get 

   \bl\label{lem-2.1jxx}   For  all $t$ sufficiently large, and $j,k\ge 1$  
 \be 
\wt Q_{j,0}(t )= \frac{ ( b t^2+t)(1-\al)  
 }{\bar t^2} \al^{j-1} \label{3.29xx}
 \ee 
and
 \be 
\wt Q_{0,k}(t)= \frac{  b t^2(1-\bb ) +\bb t
 }{\bar t^2} \bb^{j-1} \nn.\label{3.30xx}
 \ee 
Furthermore,   for all $t$ sufficiently large and for all   $1\le j\le k$\bea
\lefteqn{
\wt Q_{j,k}(t)\label{2.22js.a}}\\&& ={ \al^{j-1} \bb^{k-1}\over \bar t^2}\sum_{p=0}^{j}\wt t^{-2p}{j\choose p}{k\choose p}\nn\\
&&\quad\(-b t^2\((1-\al)(1-\bb)-{(1-\bb)p\over j}-{(1-\al)p\over k}+\wt t^{-2}{j-p\over j}{k-p\over k}
\)\right.
\nn\\
&&\qquad\qquad\left.+t \bb\(1-\al-{p\over j}  \)\)\nn.
    \eea  
\el

The analysis in Section \ref{sec-4} only uses the fact that $\ga<\ff$.  Since $b<\ff$, Lemma \ref{lem-1.1} also holds in this case.

In going through Section \ref{sec-5} we see the coefficients of $\wt Q$ change, but  they still  lead to similar inequalities that allow us to complete the proof.
 Using (\ref{3.29xx}), (\ref{3.30xx}) and (\ref{3.18}) we get 
  \be 
\wt Q_{j,0}(t )= \(\frac{b^{2}}{d^{2}t} +O(1/t)\)\al^{j-1}  \label{3.29oxx}
 \ee 
and
  \be 
\wt Q_{0,k}(t )= \(\frac{d+1}{d^{2}t} +O(1/t)\)\bb^{k-1}  \label{3.29oxxs},
 \ee 
since $ab=d+1$.

We next consider he analogue of (\ref{2.4c}) which we denote by $\wt R_{j,k}(t)$.
We see that in computing this the first two lines  of the analogue of (\ref{2.4f}) remain  unchanged, except for replacing $c$ by 1. The last two lines of  (\ref{2.4f}) are now
\bea 
   &&\(-b t^2\((1-\al)(1-\bb)-{(1-\bb)p\over j}-{(1-\al)p\over k}+\wt t^{-2}{j-p\over j}{k-p\over k}
\)\right. \nonumber\\
&&\quad\left.+t\bb\(1-\al-{p\over j} \)\).\label{7.8}
 \eea
Therefore, in place of (\ref{2.4fx}), we get     \begin{equation}
  {\wt R_{j,k,0}(t) \over  \al^{j-1}\bb^{k-1}}= {1-2c^{2}(b/d) +O\( t^{-1} \)\over d^{2} t^2}\label{2.4fxaa}.
   \end{equation}
   Using (\ref{2.4f2}), with $\ga$  replaced by $b$ and Lemma \ref{lem-2.2}, we see that (\ref{7.8})
 \begin{eqnarray}
   && =-\frac{b}{d^{2}}\((ab+1) -{ bp(dt-b) \over k}-{ap(dt-a) \over j} +{p^{2}\over jk}\)\\
   &&\qquad+\( \frac{b}{d}-{  p(dt-a)\over dj}\)+O\( {1 \over t} \) .\nn\\
      &&=-{2b\over d^{2}}	 +{   p(dt-a)\over d^{2}j}+{ b^{2} p(dt-b)\over d^{2}k}  -{b p^{2}\over jk} +O\( {1 \over t} \)\label{96ss}\\
  &&
\ge\frac{2}{d^{2}} \( p  {\sqrt{(dt-a)(dt-b)} \over \sqrt{jk} }-1\)b -{b p^{2}\over d^{2} jk} +O\( {1 \over t} \).\nn
   \end{eqnarray}
Comparing this inequality to the first line of (\ref{5.18}) we see that we have exactly what we need to complete the proof in this case.  The rest of the argument in Section \ref{sec-5} only uses the fact that $\ze>0$. It is now replaced by $b>0$. 
  Thus we get Lemma \ref{lem-5.3}   and, by Remark \ref{rem-5},    Theorem  \ref{theo-cp}, when {$  (c_{1},c_{2})=(c,0)$}} and   $EG_{1}G_{2}>0$. However this proof holds for $c$ positive or negative, so if $EG_{1}G_{2}<0$, we simply note that 
\begin{equation}
   ( G_{1} ^{2},(G_{2}+c)^{2})\stl  ( G_{1} ^{2},(-G_{2}-c)^{2}) .
   \end{equation}
Since $EG_{1}(-G_{2})>0$ we are in the case just proved so, by Remark \ref{rem-5},  Theorem  \ref{theo-cp} holds in this case also.
\qed

   \section{ Proof of Theorem  \ref{theo-cp}  for  {$\bf (c_{1},c_{2})\in R^{1}\times R^{1}$} }\label{sec-8}
   
   It is simple to complete the proof from the results already obtained. Suppose neither $c_{1}$ nor $c_{2}$ are equal to zero. Then, clearly,
   \begin{equation}
    ( (G_{1}+c c_{1}) ^{2},(G_{2}+cc_{2})^{2})
   \end{equation}
   is infinitely divisible, if and only if   
    \begin{equation}
    ( (G_{1}/c_{1}+c  ) ^{2},(G_{2}/c_{2}+c )^{2})
   \end{equation} 
      is infinitely divisible. We have already shown that there exists a critical point $\wt c>0$ such that 
     \begin{equation}
    ( (G_{1}/c_{1}+c  ) ^{2},(G_{2}/c_{2}+c )^{2})
   \end{equation} 
   is infinitely divisible for all $|c|\le \wt c$ and not for $|c|> \wt c$. Consequently $\wt c$ is also a critical point for the infinite divisibility of
   \begin{equation}
    ( (G_{1} +c_{1}c  ) ^{2},(G_{2} +c_{2}c )^{2}).
   \end{equation}    
   If $c_{1}=0$ we repeat this argument for 
   \be
       ( G_{1} ^{2},(G_{2}+cc_{2})^{2}).
   \ee
   \qed
 	 
\def\noopsort#1{} \def\printfirst#1#2{#1}
\def\singleletter#1{#1}
         \def\switchargs#1#2{#2#1}
\def\bibsameauth{\leavevmode\vrule height .1ex
         depth 0pt width 2.3em\relax\,}
\makeatletter
\renewcommand{\@biblabel}[1]{\hfill#1.}\makeatother
\newcommand{\bysame}{\leavevmode\hbox to3em{\hrulefill}\,}

 \end{document}